\documentclass{elsarticle}
\usepackage{graphicx}
\usepackage{subfigure}
\usepackage{color}
\usepackage{tikz}
\usepackage{newlfont}
\usepackage{multirow}
\usepackage{array}
\usepackage{longtable} 
\usepackage{ifthen}
\usepackage{alltt}
\usepackage{float}
\usepackage{amsmath}
\usepackage{xcolor}

\usepackage{enumerate}
 \usepackage[text={6.5in,8.5in},centering]{geometry} 
\newcolumntype{C}[1]{>{\centering\let\newline\\\arraybackslash\hspace{0pt}}m{#1}}
\newcolumntype{L}[1]{>{\raggedright\let\newline\\\arraybackslash\hspace{0pt}}m{#1}}
\newcolumntype{R}[1]{>{\raggedleft\let\newline\\\arraybackslash\hspace{0pt}}m{#1}}
\newcommand{\ba}{\begin{array} }
\newcommand{\ea}{\end{array} }
\newcommand{\bae}{\begin{eqnarray}}
\newcommand{\eae}{\end{eqnarray}}
\newcommand{\bea}{\begin{eqnarray*}}
\newcommand{\eea}{\end{eqnarray*}}
\newcommand{\be}{\begin{equation}}
\newcommand{\ee}{\end{equation}}

\newcommand{\modifyb}[1]{\textcolor{black}{#1}}
\newcommand{\modifyr}[1]{\textcolor{black}{#1}}
\newcommand{\modifyc}[1]{\textcolor{black}{#1}}
\newcommand{\modifyg}[1]{\textcolor{black}{#1}}

\newcommand{\pr}{{\bf Proof}~~}

\usepackage{amssymb}
\usepackage{amsthm}
\usepackage{ mathrsfs }
\usepackage{amsmath}
\usepackage[mathcal]{eucal}

\usepackage{natbib}
\usepackage{graphicx}
\usepackage{color}
\usepackage[colorlinks]{hyperref}
\usepackage{lscape}
\usepackage{graphicx}
\usepackage{epstopdf}
\newtheorem{theorem}{\hskip\parindent\bf Theorem}[section]

\usepackage{multicol}
\usepackage{tabularx}
\usepackage{ctable}
\usepackage{booktabs}
\usepackage{setspace}
\usepackage{bm}
\begin{document}

 \markboth{Modeling of Social Insects}{Collaborations}
\title{Disease dynamics of Honeybees with Varroa destructor as parasite and virus vector}
\author{Yun Kang\footnote{Sciences and Mathematics Faculty, College of Letters and Sciences, Arizona State University, Mesa, AZ 85212, USA ({\tt yun.kang@asu.edu})}, Krystal Blanco\footnote{Simon A. Levin Mathematical and Computational Modeling Sciences Center, Arizona State University, Tempe,  AZ 85281, USA}, Talia Davis  \footnote{School of Mathematical and Statistical Sciences, Arizona State University, Tempe, AZ 85281, USA}, Ying Wang \footnote{School of Life Sciences,
Arizona State University, Tempe, AZ 85281, USA,} and Gloria DeGrandi-Hoffman \footnote{Carl Hayden Bee Research Center,
USDA-ARS, 2000 East Allen Road, Tucson, AZ 85719.}}
\begin{abstract}
The worldwide decline in honeybee colonies during the past 50 years has often been linked to the spread of the parasitic mite \textit{Varroa destructor} and its interaction with certain honeybee viruses carried by Varroa mites. In this article, we propose a honeybee-mite-virus model that incorporates (1) parasitic interactions between honeybees and the Varroa mites; (2) five virus transmission terms between honeybees and mites at different stages of Varroa mites: from honeybees to honeybees, from adult honeybees to phoretic mites, from honeybee brood to reproductive mites, from reproductive mites to honeybee brood, and from honeybees to phoretic mites; and (3) Allee effects in the honeybee population generated by its internal organization such as division of labor.  We provide completed local and global analysis for the full system and its subsystems. Our analytical and numerical results allow us have a better understanding of the synergistic effects of parasitism and virus infections on honeybee population dynamics and its persistence. Interesting findings from our work include: (a) Due to Allee effects experienced by the honeybee population, initial conditions are essential for the survival of the colony. (b) Low adult honeybee to brood ratios have destabilizing effects on the system, generate fluctuated dynamics, and potentially lead to a \emph{catastrophic event} where both honeybees and mites suddenly become extinct. This catastrophic event could be potentially linked to Colony Collapse Disorder (CCD) of honeybee colonies. (c) Virus infections may have stabilizing effects on the system, and could make disease more persistent in the presence of parasitic mites. Our model illustrates how the synergy between the parasitic mites and virus infections consequently generates rich dynamics including multiple attractors where all species can coexist or go extinct depending on initial conditions. Our findings may provide important insights on honeybee diseases and parasites and how to best control them. 
\end{abstract}
\bigskip
\maketitle


\section{Introduction}

Honeybees are the world's most important pollinators of food crops. It is estimated that one third of food that we consume each day mainly relies on pollination by bees. For example, in the United States, honeybees are major pollinators of alfalfa, apples, broccoli, carrots and many other crops, and hence are of economic importances. Honeybees have an estimated monetary value between \$15 and \$20 billion dollars annually as commercial pollinators in the U.S \cite{Johnson:2010aa}. There are growing concerns both locally and globally that despite a 50\% growth in honeybee stocks, the supply cannot keep up with the over 300\% increase in agricultural demands \cite{Williams:2010aa}. Therefore, the recent sharp declines in honeybee populations have been considered as a global crisis. The most recent data from the 2012-2013 winter has shown an average loss of 44.8\% of hives in the U.S., and a total of 30.6\% loss of commercial hives \cite{Steinhauer:2014aa}. Some beekeepers have reported \modifyb{a loss of 90\%} of their hives \cite{Ellis:2010aa,Oldroyd:2007aa}. \\

Between 1972 and 2006, the wild honeybee populations declined severely and are now considered virtually nonexistent \cite{Moritz:2007aa,Watanabe:1994aa}. Hence the use of commercial honeybees for pollination is extremely important. Beginning in 2006, beekeepers began to report an unusual phenomenon in dying bee colonies. Worker bees would leave the colony to forage and never return, leaving the queen and the young behind to die. No dead worker bees were found at the nest sites; they simply disappear \cite{Cox-Foster:2007aa, Evans:2009aa}. This phenomenon is known as Colony Collapse Disorder (CCD), which is a serious problem threatening the health of honey bees and therefore the economic stability of commercial beekeeping and pollination operations.\\

The exact causes and triggering factors for CCD have not been completely understood yet. Researchers have proposed several possible causes of CCD including stress on nutritional diet, harsh winter conditions,  lack of genetic diversity, exposure to certain pesticides, diseases, and parasitic mites \textit{Varroa destructor} which are also vectors of viral diseases of honeybees \cite{Johnson:2010aa,Ratti:2012aa}. Even before CCD was detected in honeybee colonies, studies showed that most of the loses could be generally attributed to two main causes: the vampire mite, \text{Varroa destructor}, which feeds on host haemolymph, weakens host immunity and exposes the bees to a variety of viruses, and the tracheal mite, which infests the breathing tubes of the bee, punctures the tracheal wall and sucks the bee's blood and also exposes the bee to a variety of viruses \cite{Sammataro:2000aa,Shen:2005aa,Martin:2001aa}. Since then, Varroa mites have been implicated as the main culprit in dying colonies. For example, in Canada, Varroa mites have been found to be the main reason behind wintering losses of bee colonies \cite{Guzman-Novoa:2010aa}, and more generally studies have shown that if the mite population is not properly controlled, the honeybee colony will die  \cite{Seeley:1978aa}. Recent studies also suggest that the Varroa mite could be a contributing cause of CCD since they not only ectoparasitically feed on bees, but also vertically transmit a number of deadly viruses to the bees \cite{Ma:2009aa,Le-Conte:2010aa}. There have been at least 14 viruses found in honeybee colonies \cite{Bailey:2013aa,Ma:2009aa}, which can differ in intensity of impact, virulence, etc. for their host. For example, the Acute Bee Paralysis Virus (ABPV) affects the larvae and pupae which fail to metamorphose to adult stage, while in contrast the Deformed Wing Virus (DWV) affects larvae and pupae, which can still survive to the adult stage \cite{vanEngelsdorp:2010aa}. \\

  Mathematical models are powerful tools that could help us obtain insights on potential ecological processes that link to CCD, and important factors that contribute to the mortality of honeybees. Few sophisticated mathematical models of honeybee populations have been previously developed. DeGrandi-Hoffman \textit{et al.}   \cite{DeGrandi-Hoffman:1989aa} produced the first time-based honeybee colony growth model. Martin  \cite{Martin:1998aa} developed a simulation model consisting of ten components, which linked together various aspects of mite biology using computer software (ModelMaker); and Martin  \cite{Martin:2001aa}  later extended this model by including a bee model adapted from \cite{DeGrandi-Hoffman:1989aa} to explain the link between the Varroa mite and collapse of the host bee colony. Wilkinson and Smith  \cite{Wilkinson:2002aa} proposed a difference equation model of Varroa mites reproducing in a honeybee colony. Their study focused on parameter estimations and sensitivity analysis. Simulation models are useful but may be too complex to study mathematically and obtain general predictions. \\ 

More recently, mathematical models have been formulated to explore potential mechanisms causing CCD to the honeybee colony. \modifyg{it would be helpful to include the VARROAPOP model here because you can later compare the behavior and trends from your model with previously published models.}  Sumpter and Martin  \cite{Sumpter:2004aa} modeled the effects of a constant population of Varroa mites on the brood and on adult worker bees, and found that sufficiently large mite infestations may make hives vulnerable to collapse from viral epidemics. Eberl \textit{et al.}  \cite{Eberl:2010aa} developed a model connecting Varroa mites to CCD by including brood maintenance terms which reflect that a certain number of worker bees is always required to care for the brood in order for them to survive. They found an important threshold for the number of hive worker bees needed to maintain and take care of the brood.  Khoury \textit{et al.} \cite{Khoury:2011aa,Khoury:2013aa} developed differential equations models to study different death rates of foragers and the impact it had on colony growth and development. They then linked their results to CCD. Betti \textit{et al.} \cite{Betti:2014aa} studied a model that combines the dynamics of the spread of disease within a bee colony with the underlying demographic dynamics of the colony to determine the ultimate fate of the colony under different scenarios. Their results suggest that the age of recruitment of hive bees to foraging duties is a good early marker for the survival or collapse of a honeybee colony in the face of infection. Kribs-Zaleta \textit{et al.}  \cite{Kribs-Zaleta:2014aa} created a model to account for both healthy hive dynamics and hive extinction due to CCD, modeling CCD via a transmissible infection brought to the hive by foragers. Perry \textit{et al.} \cite{Perry:2015aa} examined the social dynamics underlying the dramatic colony failure with an aid of a honeybee population model. Their model includes bee foraging performance varying with age, and displays dynamics of colony population collapse that are similar to field reports of CCD. These models, no doubt, are insightful and provide us a better understanding on the potential mechanisms that link to CCD. However, most of these models only account for the honeybee population dynamics with mites or viruses but not both. \\ 

The host-parasite relationship between honeybees and Varroa mites has been complicated by the mite's close association with a  wide range of honeybee viral pathogens. In order to understand how Varroa mite infestations and the related viruses transmitted to honeybees affect honeybee population dynamics, and which may link to CCD, there is a need to develop realistic and mathematically tractable models that include both mite and pathogen population dynamics.  The goal of our work is to develop a useful honeybee-mite-virus system to obtain better understanding on the synergistic effects of honeybee-mite interactions and honeybee-virus interactions on the honeybee populations dynamics, thus develop good practices to control these parasites to maintain or increase honeybee population. The most relevant modeling papers for our study purposes are by Sumpter and Martin \cite{Sumpter:2004aa}, and Ratti \textit{et al.} \cite{Ratti:2012aa} whose work examined the transmission of viruses via Varroa mites, using the susceptible-infectious (SI) disease modeling framework with mites as vectors for transmission. However, Sumpter and Martin assumed that the mites' population is constant while Ratti \textit{et al.} took no account of the fact that virus transmissions occur at different biological stages of Varroa mites and honeybees.\\

 In this article, we follow both approaches of Sumpter and Martin \cite{Sumpter:2004aa} and Ratti \textit{et al.} \cite{Ratti:2012aa}, and propose a honeybee-mite-virus model that incorporates (1) parasitic interactions between honeybee and Varroa mites; (2) different virus transmission terms that account for the virus transmission among honeybees, between honeybees and mites at different stages of Varroa mites; and (3) Allee effects in the honeybee population generated by the internal organization of honeybees, including division of labor. Our proposed model will allow us explore the following questions:\\
\begin{enumerate}
\item What are the dynamics of a system only consisting of honeybees and the disease?
\item What are the dynamics of a system only consisting of honeybees and Varroa mites?
\item What are the synergistic effects of Varroa mites and the disease on the honeybee population, and how may these synergistic effects contribute to CCD?
\item  How can we maintain honeybee populations?\\
\end{enumerate}

The structure of the \modifyb{remainder of the} article is organized as follows: In Section 2, we first provide the biological background of honeybees, Varroa mites, and the associated virus transmission routes in the honeybee-mite system; then we derive our SI-type model for honeybees co-infected with the mite and virus. In Section 3,  we perform local and global analysis of the proposed model and the related subsystems. The results from the analysis are then connected to biological contexts and implications. Additionally, we also explore numerical simulations of the subsystems and the full system to obtain the effects of each parameter in our system. In Section 4, we summarize our results and the related biological implications of our studies in finding potential causes of Colony Collapse Disorder.  We also provide potential projects for future work. The detailed mathematical proofs of our theoretical results are provided in the last section.\\


\section{Biological background and model derivations}\label{sec:biological_background}

\textbf{Honeybee colony:} During the spring and summer, a honeybee colony typically consists of a single reproductive queen, 20,000 -- 60,000 adult worker bees, 10,000 -- 30,000 individuals at the brood stage (egg, larvae and pupae) and up to hundreds of male drones. During the winter, the colony typically reduces in size and consists of a single queen and somewhere between 8,000 -- 15,000 worker bees  \cite{Martin:2001aa}. A large population of workers carry out the tasks of the bee colony, which include foraging, pollination, honey production and, in particular, caring for the brood and rearing the next generation of bees. The queen is the only fertile individual of the colony and has an average life span of 2 -- 3 years \cite{Sumpter:2004aa}. During the peak season (in the summer), the queen lays up to 2000 eggs per day, where fertilized eggs produce female worker bees, or much more rarely queens, while drones develop from non-fertilized eggs \cite{Betti:2014aa}. The bees go through the following stages in development: egg (about 3 days), larvae (about 7 days), pupae (about 14 days), and adult.  The life span of an adult worker bee also depends on the season. Workers usually have a lifespan of 3 -- 6 weeks during the spring and summer, and are reported to live as long as 4 months during the winter \cite{Page:2001aa}. The adult drone life span is  typically 20 -- 40 days, with reports of drone living up to 59 days under optimal colony conditions \cite{Page:2001aa,Howell:1933aa}.\\

Let $N_h(t)$ be the total number of honeybees in the colony, including the larvae, pupae and adult bee (both hives and foragers) at time $t$. \modifyb{The subscript $h$ means honeybees for all future notations}. \modifyc{Define $\xi_h\in [0,1]$ as the percentage of adult honeybee population in the colony, then $(1-\xi_h) N_h$ is the brood population and $\frac{\xi_h}{1-\xi_h}$ is the ratio of adult honeybees to the brood in colony. Empirical study shows that the successful honeybee colony should have $\frac{\xi_h}{1-\xi_h}>2$ (see \cite{Schmickl:2007aa}). We should expect that the value of $\xi_h$ varies with time. In our current model, instead of employing the explicit age structure model, we let $\xi_h$ be a parameter. Though brood and adult numbers change throughout the year, the ratio remains fairly constant and is a limiting factor in proportion of eggs that are reared into larvae and emerge as adults. An addition justification of this simplification is that  the ratio of adult bees $\xi_h$ is indeed a constant at the steady state of the explicit age structure model. Under this simplification, we are able to obtain analytical results on how $\xi_h$ affects dynamics of the model with essential biological components.}
In the absence of mites and virus, the population dynamics of $N_h(t)$ is described by the following nonlinear equation:
\bae\label{Hpopulation}
N_h'&=& \frac{r (\xi_hN_h)^2}{K+(\xi_hN_h)^2}-d_h N_h\eae where \modifyb{$'$ is the sign of the derivative with respect time}; the parameter $r$ is the maximum birth rate, specified as the number of worker bees born per day; the parameter $\sqrt{K}$ is the size of the bee colony at which the birth rate is half of the maximum possible rate; and $d_h$ is the average death rate of the worker honeybees. The term $\frac{(\xi_hN_h)^2}{K+(\xi_hN_h)^2}$ describes that the successful survival of an egg which will develop into a worker bee needs the care of adult honeybees ($\xi_hN_h$) inside the colony and also needs food brought in by the honeybee foragers. \modifyr{This term also implicitly assumes that more population of adult honeybees inside the colony can increase the survival of an egg developing into a worker which is supported by the empirical study showed in \cite{Schmickl:2007aa}.} Our modeling  approach of \eqref{Hpopulation} follows the modeling idea in \cite{Eberl:2010aa} for honeybee diseases and in \cite{Kang:2011aa} for the population of leaf-cutter ants. This term implicitly includes the internal organization of the honeybee population, such as division of labor. \modifyr{Model \eqref{Hpopulation} implies that honeybee population is able to persist when $\xi_h$ is above an threshold, i.e., $\xi_h>\frac{2d_h\sqrt{K}}{r}$. More detailed analysis and related results are provided in the next section. }\\

\textbf{Varroa mites:} Varroa mites were \modifyb{first reported in Kentucky in the Bluegrass region of the Commonwealth in 1991 in U.S}. They have since spread to become \modifyb{a major pest of honeybees in many states of U.S.} \cite{Bessin:2013aa}. Varroa mites are external honeybee parasites that attack both adult honeybees and brood, with a distinct preference for drone brood \cite{Peng:1987aa}. They suck the blood from both the adults and the developing brood, weakening them and shortening the life span of the bees which they feed on. Emerging brood may be born with deformed wings. Untreated infestations of Varroa mites can cause honeybee colonies to collapse \cite{Martin:1994aa}. \\

The mites go through a series of stages: larva, protonymph, deutonymph and then adult. Adult females undergo two phases in their life cycle, the phoretic and reproductive phases. During the phoretic phase, female Varroa feed on adult bees and are passed from bee to bee as they pass one another in the colony. During the phoretic phase, the female Varroa mites live on adult bees and can usually be found between the abdominal segments of the bees. The mites puncture the soft tissue between the segments and feed on bee hemolymph, harming the host \cite{Rosenkranz:2010aa,Calderon:2007aa}. Mite reproduction can occur only if honeybee brood is available. A female mite enters the brood cell about one day before capping and will be sealed in with the larva. After the capping of the cell, it lays a single male egg and several female eggs at 30-hour intervals \cite{Vetharaniam:2006aa}, and the mite feeds and develops on the maturing bee larva. When the host bee leaves the cell, the mature female mites leave the cell. The male mite dies after mating with his sisters, and if immature female mites are present they die as they come out of the cell, as they cannot survive once outside the cell. The adult female mite begins searching for other bees or larvae to parasitize.\\

The phoretic period of the mite appears to contribute to the mite's reproductive ability, which may last 4.5 to 11 days when brood is present in the hive; or as long as five to six months during the winter when \modifyb{little or no brood} is present in the hive. Consequently, female mites living when brood is present in the colony have an average life expectancy of 27 days, yet in the absence of brood, they may live for many months. In the average temperate climate, mite populations can increase 12-fold in colonies which have brood half of the year and 800-fold in colonies which have brood year-round. \modifyb{This period usually begins in late winter when brood rearing resumes from a winter period when little or no brood is present. The period of mite population increase continues through the spring and summer and peaks in the fall when brood rearing is nearly done.} This makes the mites very difficult to control, especially in warmer climates where colonies maintain brood year-round \cite{Ellis:2010ab}.\\

Let $N_{m} (t)$ be the number of adult female Varroa mites in the honeybee colony in the absence of  a virus where \modifyb{the subscript $m$ means Varroa mites for all future notations}.  Varroa mites feed on the haemolymph of brood and adult honeybees, and their reproduction depends on the availability of the brood and the population of the reproductive female Varroa mites. \modifyc{Similarly to the case of honeybees, we incorporate an implicit age structure model of Varroa mites by defining a parameter $\xi_m \in [0,1]$ as the percentage of Varroa mites at the  phoretic stage. This implies that the honeybee colony has the population size of $(1-\xi_m) N_m$ as the reproductive Varroa mites and the population size of $\xi_m N_m$ as the phoretic Varroa mites. This simplification still allows us to investigate the parasitic interactions between Varroa mites and honeybees rigorously with essential biological components.}\\

We model the parasitic interactions between Varroa mites and honeybees by using the Holling Type I functional responses, i.e. \modifyb{$\hat{\alpha}(1-\xi_h) N_h(1-\xi_m) N_m=\alpha N_h N_m$ where $\hat{\alpha}$ is the parasitism rate; the term $(1-\xi_h) N_h$ is the honeybee brood population; the term $(1-\xi_m) N_m$ is the reproductive female Varroa mites population; and $\alpha=\hat{\alpha} (1-\xi_h)(1-\xi_m)$.} \modifyr{Therefore, in the presence of Varroa mites $N_m$, the dynamics of the honeybee population $N_h$ can be described as:
$$\begin{array}{lcl}
			N_h'&=&  \frac{r (\xi_hN_h)^2}{K+(\xi_hN_h)^2}-d_h N_h-\alpha N_h N_m=\frac{r (\xi_hN_h)^2}{K+(\xi_hN_h)^2}-(d_h +\alpha N_m)N_h
			\end{array}$$which implies that the parasitism decreases the life span of the honeybee, i.e., the average life span has been reduced to $\frac{1}{d_h +\mathbf{\alpha N_m}}$ after parasitism. Here we do not assume that parasitism would lead to the death of the honeybees for sure. The Varroa mite population depends on the nutrient obtained from honeybees, thus, the dynamics of Varroa mites population can be described as:
			$$\begin{array}{lcl}
			N_m'&=&c \alpha N_hN_m-d_m N_m\end{array}$$ where the parameter $\alpha$  measures the parasitic rate of Varroa mites; $c$ is the conversion rate from  nutrient consumption obtained from honeybees to sustenance for Varroa mites reproduction; and $d_m$ is the natural death rate of Varroa mites. Therefore, in the absence of virus, the population dynamics of Varroa mites and  honeybees  can be described by the following two nonlinear equations:
            \bae\label{SMpopulation}
		\begin{array}{lcl}
			N_h'&=&  \frac{r (\xi_hN_h)^2}{K+(\xi_hN_h)^2}-d_h N_h-\alpha N_h N_m\\\\
			N_m'&=&c \alpha N_hN_m-d_m N_m
			\end{array}.
			\eae The modeling approach of Model \eqref{SMpopulation} has applied the traditional host-parasite modeling framework including non-lethal parasites \cite{Anderson1978, AndersonMay1981}.}
Model \eqref{SMpopulation} also implies that Varroa mites population $N_m$ goes extinct if the population of honeybees $N_h$ goes extinct. \\

\textbf{Varroa mites act as a disease-vector for virus transmissions:} 
Varroa mites not only feed on host haemolymph and weaken host immunity, but they also expose honeybee colonies to at least 14 different viruses including deadly viruses such as ABPV and DWV. Varroa mites can transmit viruses in their reproduction phase to honeybee brood and during the phoretic phase to adult honeybees. To model the virus transmission between Varroa mites and honeybees during these two phases, we let $S_h(t), S_m(t)$ be the susceptible population of honeybees and Varroa mites, respectively; and $I_h(t), I_m(t)$ be the \modifyb{virus infected} population of honeybees and Varroa mites, respectively. Then the total population of honeybees is $N_h(t)=S_h(t)+I_h(t)$, and the total population of Varroa mites is $N_m(t)=S_m(t)+I_m(t)$.\\

The virus transmission between female Varroa mites and honeybees can occur in the following two phases of the Varroa mite life cycle:
\begin{enumerate}
\item  The honeybee colony has $\xi_h S_h$ susceptible adult honeybees;  $\xi_h I_h$ \modifyb{virus infected} adult honeybees, $\xi_m S_m$ susceptible phoretic female Varroa mites; and $\xi_m I_m$ \modifyb{virus infected} phoretic female Varroa mites. In the phoretic phase, female Varroa mites move between adult bees both spontaneously and just prior to the death of their host bee \cite{Sumpter:2004aa}. Following the approach of  \cite{Martin:2001aa}, we assume that virus transmission is frequency dependent, i.e., 
\begin{itemize}
\item \modifyr{We model the rate at which susceptible adult honeybees are \modifyb{virus infected} by the infected phoretic female Varroa mites (IPFM) based on the approach of \cite{Martin:2001aa, Sumpter:2004aa, Eberl:2010aa, Ratti:2012aa}. This rate can be described as follows:\\
$$\begin{array}{ll}
&\underbrace{\beta_{mh}}_{\text{probability being infected after contacts}}\times \underbrace{\xi_m I_m}_{\text{population of infected IPFM}} \times \underbrace{\frac{\xi_h S_h }{\xi_h S_h+\xi_h I_h}}_{\text{probability contacts susceptible adult honeybees}}\\\\
&=\frac{\beta_{mh} \xi_m S_h I_m}{S_h+I_h}.
\end{array}$$ which also implies that susceptible honey bees become virus infected at a rate proportional to the ratio of the population of the infected phoretic mites to the  population of the total honeybees. }
\item The rate at which susceptible phoretic female Varroa mites (SPFM) are \modifyb{virus infected} by the \modifyb{virus infected} adult honeybees (IAH) is:\\
$$\begin{array}{ll}
&\underbrace{\beta_{hm}}_{\text{probability being infected after contacts}}\times \underbrace{\xi_m S_m}_{\text{population of SPFM}} \times \underbrace{\frac{\xi_h I_h }{\xi_h S_h+\xi_h I_h}}_{\text{probability being contacted by IAH}}\\\\
&=\frac{\beta_{hm} \xi_m S_m I_h}{S_h+I_h}.
\end{array}$$
\end{itemize}
\item The honeybee colony has $(1-\xi_h) S_h$ susceptible honeybee brood;  $(1-\xi_h) I_h$ \modifyb{virus infected} honeybee brood, $(1-\xi_m) S_m$ susceptible reproductive female Varroa mites; and $(1-\xi_m) I_m$ \modifyb{virus infected}  reproductive female Varroa mites. \modifyr{We do not assume that honeybee brood will die from parasitism, thus, parasitised honeybee brood or newborn Varroa mites will face virus infection if either honeybee brood or the reproductive female Varroa mite is virus infected}.
 Chen \textit{et al.} found that there is a direct relationship between virus frequency and the number of mites to which honeybee brood were exposed, i.e., the more donor mites that were introduced per cell, the greater the incidence of virus that was detected in the honeybee brood \cite{Chen:2006aa,Chen:2006ab}. This implies that the virus transmission rate between Varroa mites and the honeybee brood during the reproductive phase of mites is density dependent, i.e., similar to the term that describes the parasitic interaction between mites and honeybee. Therefore, we have as follows:
\begin{itemize}
\item A newborn honeybee becomes \modifyb{virus-infected if it is parasitized} by the infected reproductive female Varroa mites. Thus, the rate at which susceptible honeybee brood is \modifyb{virus-infected} by the \modifyb{virus infected} reproductive female Varroa mites (IRFM) is:\\
$$\begin{array}{ll}
&\underbrace{\beta_{mh2}}_{\text{probability being infected after contacts}}\times \underbrace{(1-\xi_h) S_h}_{\text{population of healthy honeybee brood}} \times \underbrace{\hat{\alpha} (1-\xi_m)I_m}_{\text{parasitism by IRFM}}\\\\
&=\beta_{mh2}\alpha S_hI_m.
\end{array}$$

\item The reproduction of Varroa mites depends on honeybee brood. The newborn Varroa mites become \modifyb{virus infected} if either the brood  or the female Varroa mites is \modifyb{virus infected}. Thus, based on the formulation of the host-parasite interaction model \eqref{SMpopulation}, the rate at which infected newborn female Varroa mites (INFM) become \modifyb{virus infected} depends on the parasitic interaction between mites and honeybees which can be described as \modifyb{$c\alpha\left[I_h (S_m+I_m)+S_hI_m\right]$ where the term $I_h (S_m+I_m)$ is the newborn Varroa mites infected with virus through virus infected honeybee brood; and the term $S_hI_m$ is the newborn Varroa mites infected with virus through virus infected reproductive Varroa mites.}

\end{itemize}
\end{enumerate}

\textbf{The virus transmission among honeybees:} The proportion of honeybees which can infect themselves is also dependent on the total number of susceptible and \modifyb{virus infected} bees present in the colony, and hence frequency-dependent transmission is used \cite{Martin:2001aa}, which is described as follows:

$$\begin{array}{ll}
&\underbrace{\beta_{h}}_{\text{probability being infected after contacts}}\times \underbrace{ S_h}_{\text{the healthy honeybee population}} \\
&\times \underbrace{\frac{ I_h }{ S_h+ I_h}}_{\text{probability of contacting or being contacted by infected honeybees}}=\frac{\beta_{h} S_h I_h}{S_h+ I_h}.\end{array}$$\\

\textbf{The reduced fitness of honeybees due to virus infections:} The parasitic Varroa mites have been shown to act as a vector for a number of viruses including DWV, ABPV, Chronic Bee Paralysis Virus (CBPV), Slow Bee Paralysis Virus (SPV), Black Queen Cell Virus (BQCV), Kashmir Bee Virus (KBV), Cloudy Wing Virus (CWV), and Sacbrood Virus (SBV) \cite{Martin:1998aa,Martin:2001aa,Martin:2012aa,Le-Conte:2010aa}.
These virus infections  contribute to morphological deformities of honeybees such as small body size, shortened abdomen and deformed wings, which reduce vigor and longevity, and they can also influence flight duration and the homing ability of foragers \cite{Le-Conte:2010aa}.  In our model, we assume that the \modifyb{virus infected} adult honeybee population $\xi_h I_h$ contributes to the reproduction of honeybees with a reduced rate $\rho \in (0,1)$, therefore, the healthy honeybee population $S_h$ can be modeled as follows,

{\small\textbf{\bae
		\begin{array}{lcl}
			S_h'&=& \underbrace{\frac{r\xi_h^2\left(S_h+\rho I_h\right)^2}{K+\xi_h^2\left(S_h+\rho I_h\right)^2}}_{\text{ reproduction of honeybees}}-\underbrace{\frac{\beta_{h} S_hI_h}{S_h+I_h}}_{\text{honey bee infected by themselves}}-\underbrace{\alpha S_h (S_m+I_m)}_{\text{parasitism by mites}}\\\\
			&&-\underbrace{\modifyr{\frac{\beta_{mh}(\xi_m S_h)I_m}{S_h+I_h}}}_{\text{adult honeybees infected  by the phoretic mites}}
			-\underbrace{\beta_{mh2}\alpha S_hI_m}_{\text{honeybee brood infected by the reproductive mites}}-d_h S_h
		\end{array}.
		\eae} }

\noindent And the virus infected honeybee population can be modeled by the following equation,

{\small\textbf{\bae
		\begin{array}{lcl}
			I_h'=S_h\left[\frac{\beta_{h}I_h}{S_h+I_h}+\modifyr{\frac{\beta_{mh}\xi_m I_m}{S_h+I_h}}+\beta_{mh2}\alpha I_m\right]-\underbrace{\alpha I_h (S_m+I_m)}_{\text{Consumed by mites}}\\\\
			-(\underbrace{d_h}_{\text{natural honeybee mortality rate}}+\underbrace{\mu_h}_{\text{additional death due to virus infections}})I_h
		\end{array}.
		\eae} }

Let $\mu_m$ be the additional death rate of Varroa mites due to virus infections. Then the population of healthy Varroa mites $S_m$ and the \modifyb{virus infected} Varroa mites $I_m$ can be described by the following set of nonlinear equations:
{\small\textbf{\bae
		\begin{array}{lcl}
			S_m'&=&S_m\left[c\alpha S_h-\underbrace{\frac{\beta_{hm}\xi_mI_h}{S_h+I_h}}_{\text{the phoretic mites infected by adult honeybees}}-\underbrace{d_m}_{\text{natural mortality rate of mites}}\right]\\\\
				I_m'&=&c\alpha \underbrace{\left[I_h (S_m+I_m)+S_hI_m\right]}_{\text{mites born with virus infections}}+\frac{\beta_{hm}I_h(\xi_mS_m)}{S_h+I_h}-(d_m +\underbrace{\mu_m}_{\text{additional death due to virus infections}})I_m
		\end{array}.
		\eae} }

\noindent Based on the discussions above, the full model of honeybee-mites-virus population dynamics is therefore modeled by the following system of differential equations:
{\small\textbf{\bae\label{HBD}
\begin{array}{lcl}
S_h'&=& \frac{r\xi_h^2\left(S_h+\rho I_h\right)^2}{K+\xi_h^2\left(S_h+\rho I_h\right)^2}-d_h S_h-\frac{\beta_{h} S_hI_h}{S_h+I_h}-\modifyr{\frac{\beta_{mh}(\xi_m S_h)I_m}{S_h+I_h}}\\\\
&&-\beta_{mh2}\alpha I_m-\alpha S_h (S_m+I_m)\\\\
I_h'&=&S_h\left[\frac{\beta_{h}I_h}{S_h+I_h}+\modifyr{\frac{\beta_{mh}\xi_m I_m}{S_h+I_h}}+\beta_{mh2}\alpha I_m\right]-\alpha I_h (S_m+I_m)-(d_h+\mu_h)I_h\\\\
S_m'&=&S_m\left[c\alpha S_h-\frac{\beta_{hm}\xi_mI_h}{S_h+I_h}-d_m\right]\\\\
I_m'&=&c\alpha \left[I_h (S_m+I_m)+S_hI_m\right]+\frac{\beta_{hm}I_h(\xi_mS_m)}{S_h+I_h}-(d_m +\mu_m)I_m
\end{array}.
\eae} } 

For convenience, let $\hat{K}=\frac{K}{\xi_h^2}, \,\,\hat{\beta}_{mh}=\beta_{mh}\xi_m, \,\,\tilde{\beta}_{mh}=\beta_{mh2}\alpha, \,\,\hat{\beta}_{hm}=\beta_{hm}\xi_m.$ Then the full model \eqref{HBD} can be rewritten as the following model 
\bae\label{newHBD}
\begin{array}{lcl}
S_h'&=& \frac{r\left(S_h+\rho I_h\right)^2}{\hat{K}+\left(S_h+\rho I_h\right)^2}-d_h S_h-\frac{\beta_{h} S_hI_h}{S_h+I_h}-\modifyr{\frac{\hat{\beta}_{mh} S_hI_m}{S_h+I_h}}-\tilde{\beta}_{mh}S_hI_m-\alpha S_h (S_m+I_m)\\\\
I_h'&=&S_h\left[\frac{\beta_{h}I_h}{S_h+I_h}+\modifyr{\frac{\hat{\beta}_{mh}I_m}{S_h+I_h}}+\tilde{\beta}_{mh}I_m\right]-\alpha I_h (S_m+I_m)-(d_h+\mu_h)I_h\\\\
S_m'&=&S_m\left[c\alpha S_h-\frac{\hat{\beta}_{hm}I_h}{S_h+I_h}-d_m\right]\\\\
I_m'&=&c\alpha\left[I_h (S_m+I_m)+S_hI_m\right]+\frac{\hat{\beta}_{hm}I_hS_m}{S_h+I_h}-(d_m +\mu_m)I_m
\end{array}
\eae where $\alpha>0, c>0$, $\rho\in [0,1]$ and the virus transmission rates $\beta_{h},\hat{\beta}_{mh}, \hat{\beta}_{hm}, \tilde{\beta}_{mh}\in (0,1)$.  In summary, the full honeybee-mite-virus model \eqref{newHBD} incorporates (1) Allee effects of honeybees due to the cooperation of the internal organization; (2) parasitism interactions between honeybee and mites; (3) the vertical disease transmission mode modeled by the frequency-dependent disease transmission function during Varroa mites' phoretic phase; (4) the horizontal disease transmission mode modeled by the density-dependent disease transmission function during Varroa mites' reproductive phase; and (5) the reduced fitness of honeybees due to virus infections.  \modifyc{The model \eqref{newHBD} allows us to investigate the following scenarios:\begin{enumerate}
\item In the absence of the Varroa mites and virus, how population of honeybees  may persist.
\item In the absence of the virus, how the Varroa mites may affect the population dynamics of honeybees.
\item In the absence of the Varroa mites, how virus infections may affect the population dynamics of honeybees. This case can apply to the situations that honeybees are infected by virus through ecological processes such as foraging other than parasitism by Varroa mites.
\item In the presence of both the Varroa mites and virus infections, which conditions can lead to the extinction of Varroa mites, virus infections, and honeybees; and which conditions can guarantee the persistence of honeybee population.
\end{enumerate}The rest of this manuscript will focus on the dynamics of Model \eqref{newHBD} and the related subsystems.}\\


\section{Mathematical analysis}
Let $\frac{S_hI_h}{S_h+I_h}\big\vert_{S_h=I_h=0}=0$ and $\frac{I_m}{S_m+I_m}\big\vert_{S_m=I_m=0}=0$. Define $X=\{(S_h, I_h, S_m, I_m)\in\mathbb R^4_+: S_h+ I_h>0 \mbox{ and } S_m+ I_m>0\}$, then $X$ can be considered as the state space of our model \eqref{newHBD}. To continue the analysis, \modifyc{let us define $N_h=S_h+I_h, N_m=S_m+I_m$ and $N= cN_h+N_m$ as the population of honeybees, the population of mites, and the sum of the population of honeybees and mites, respectively. In addition, we let $d=\min\{d_h,d_m\}$, and define $N^*$ as the upper bound of the  sum of the population of honeybees and mites and $N^c$ as the corresponding threshold, where $$\begin{array}{lll}
N^c=c\frac{\frac{r}{d}-\sqrt{\left(\frac{r}{d}\right)^2-4\hat{K}}}{2},&N^*=c\frac{\frac{r}{d}+\sqrt{\left(\frac{r}{d}\right)^2-4\hat{K}}}{2}\end{array}.$$
We let $\bar{N}^*_h, \underline{N}^*_h$ be the upper bound, lower bound of the population of honeybees, respectively, and $\bar{N}^c_h, \underline{N}^c_h$ be the corresponding thresholds, where
$$\begin{array}{lll}
\bar{N}^c_h=\frac{\frac{r}{d_h}-\sqrt{\left(\frac{r}{d_h}\right)^2-4\hat{K}}}{2},&\bar{N}^*_h=\frac{\frac{r}{d_h}+\sqrt{\left(\frac{r}{d_h}\right)^2-4\hat{K}}}{2}\\\\
\underline{N}^c_h= \frac{\frac{r}{d_h +\mu_h+\alpha N^*}-\sqrt{\left(\frac{r}{d_h +\mu_h+\alpha N^*}\right)^2-4\hat{K}/\rho^2}}{2},&\underline{N}^*_h= \frac{\frac{r}{d_h +\mu_h+\alpha N^*}+\sqrt{\left(\frac{r}{d_h +\mu_h+\alpha N^*}\right)^2-4\hat{K}/\rho^2}}{2}
\end{array}.$$
And we let $S_h^*$ be the lower bound of the population of susceptible honeybees, and $S_h^c$ be the corresponding threshold, where
$$\begin{array}{lll}
S_h^c=\frac{\frac{r}{d_h+\beta_h+\modifyr{\frac{\hat{\beta}_{mh} N^*}{\underline{N}^*_h}}+(\tilde{\beta}_{mh}+\alpha)(N^*-c\underline{N}^*_h)}-\sqrt{\left(\frac{r}{d_h+\beta_h+\modifyr{\frac{\hat{\beta}_{mh} N^*}{\underline{N}^*_h}}+(\tilde{\beta}_{mh}+\alpha)(N^*-c\underline{N}^*_h)}\right)^2-4\hat{K}}}{2},\\\\
S_h^*=\frac{\frac{r}{d_h+\beta_h+\modifyr{\frac{\hat{\beta}_{mh} N^*}{\underline{N}^*_h}}+(\tilde{\beta}_{mh}+\alpha)(N^*-c\underline{N}^*_h)}+\sqrt{\left(\frac{r}{d_h+\beta_h+\modifyr{\frac{\hat{\beta}_{mh} N^*}{\underline{N}^*_h}}+(\tilde{\beta}_{mh}+\alpha)(N^*-c\underline{N}^*_h)}\right)^2-4\hat{K}}}{2}
\end{array}.$$
Define $f^b(x,y)=\frac{\frac{r}{x}+\sqrt{\left(\frac{r}{x}\right)^2-4\hat{K}/y}}{2}$ and $f_b(x,y)=\frac{\frac{r}{x}-\sqrt{\left(\frac{r}{x}\right)^2-4\hat{K}/y}}{2}$, then we have
$$\frac{\partial f^b(x,y)}{\partial x}<0,\,\frac{\partial f^b(x,y)}{\partial y}<0,\,\frac{\partial f_b(x,y)}{\partial x}>0,\,\frac{\partial f_b(x,y)}{\partial y}>0$$ which imply the following inequalities
$$S_h^c<{S}_h^c<S_h^*<\bar{N}^*_h,\,\, N^c\leq \bar{N}^c_h<\underline{N}^c_h< \underline{N}^*_h<\bar{N}_h^*\leq N^*/c.$$} 

\begin{theorem}[Basic dynamical properties]\label{th1:basics} Assume that all parameters are strictly positive and $\rho\in[0,1]$.  The model \eqref{newHBD} is positively invariant and bounded in the state space $X$, which is attracted to the following compact set 
$$C=\{(S_h, I_h, S_m, I_m)\in\mathbb R^4_+: 0\leq c(S_h+ I_h)+ (S_m+ I_m)=cN_h+N_m\leq N^*\}$$ provided that $\frac{r}{d}>2\sqrt{\hat{K}}$ and time is large enough. Moreover, the following statements hold for Model \eqref{newHBD}:
\begin{itemize}
\item If $ \frac{r}{2\sqrt{\hat{K}}}>d_h$, then the total population of honeybees $N_h$ is bounded by $\bar{N}_h^*$, i.e., 
$$\limsup_{t\rightarrow\infty} N_h(t)\leq \bar{N}^*_h.$$ If $\frac{r}{2\sqrt{\hat{K}}}>\frac{d_h +\mu_h+\alpha N^*}{\rho}\mbox{ and } N_h(0)>\underline{N}^c_h$ hold, then the total population of honeybees $N_h$ is persistent, i.e., 
$$\underline{N}^*_h\leq\liminf_{t\rightarrow\infty} N_h(t)\leq \limsup_{t\rightarrow\infty} N_h(t)\leq \bar{N}^*_h\leq\limsup_{t\rightarrow\infty} N(t)/c=N^*/c.$$
\item 
If the inequalities  $\frac{r}{2\sqrt{\hat{K}}}>\max\Big\{d_h+\beta_h+\modifyr{\frac{\hat{\beta}_{mh} N^*}{\underline{N}^*_h}}+(\tilde{\beta}_{mh}+\alpha)(N^*-c\underline{N}^*_h),\, \frac{d_h +\mu_h+\alpha N^*}{\rho}\Big\}$ \mbox{ with } $ N_h(0)\geq S_h(0)>S_h^c $ hold, 
then $S_h$ is persistent with the following properties:
$$S^*_h\leq \liminf_{t\rightarrow\infty} S_h(t)\leq \liminf_{t\rightarrow\infty} N_h(t)\leq \liminf_{t\rightarrow\infty} N(t)/c\leq \limsup_{t\rightarrow\infty} N(t)/c\leq N^*/c.$$
\item The extinction equilibrium $E_0=(0,0,0,0)$ is always local stable. Moreover, the system \eqref{newHBD} converges to $E_0$ globally if $d_h>\frac{r}{2\sqrt{\hat{K}}}$ holds; and it converges to $E_0$ locally if the initial population satisfies either $N(0)<N^c$ or $N_h(0)<\bar{N}^c_h$.\\
\end{itemize}
\end{theorem}

\noindent\textbf{Notes:} The positive invariance and boundedness results from Theorem \ref{th1:basics} imply that our model is well-defined biologically. In addition, Theorem \ref{th1:basics} indicate follows:
\begin{enumerate}
\item Initial conditions are important for the persistence of honeybees.
\item The inequality $\frac{r}{2\sqrt{\hat{K}}}>d_h$ is a necessary condition for honeybee persistence, i.e., the large intrinsic growth rate $r$, small half saturation $\hat{K}$, and the small death rate of honeybees $d_h$.
\item The small values of disease transmission rates $\beta_h$, $\hat{\beta}_{mh}$, $\tilde{\beta}_{mh}$; and small values of mite attacking rate $\alpha$ are also important for the persistence of the healthy honeybee population $S_h$.\\
\end{enumerate}
\modifyr{Recall that $d=\min\{d_h,d_m\}, N^c=c\frac{\frac{r}{d}-\sqrt{\left(\frac{r}{d}\right)^2-4\hat{K}}}{2},\,N^*=c\frac{\frac{r}{d}+\sqrt{\left(\frac{r}{d}\right)^2-4\hat{K}}}{2}$ and\vspace{-8pt}
$$\underline{N}^*_h= \frac{\frac{r}{d_h +\mu_h+\alpha N^*}+\sqrt{\left(\frac{r}{d_h +\mu_h+\alpha N^*}\right)^2-4\hat{K}/\rho^2}}{2}.$$
Theorem \ref{th1:basics} implies that under proper initial conditions, honeybees can persist if $\frac{r}{2\sqrt{\hat{K}}}>\frac{d_h +\mu_h+\alpha N^*}{\rho}\mbox{ and } N_h(0)>\underline{N}^c_h$ holds. Notice that $\xi_h$ is the ratio of adult honeybees in the colony, and one of sufficient conditions that guarantee the persistence of honeybee population is the following inequality:
$$\frac{r}{2\sqrt{\hat{K}}}=\frac{r{\xi_h}}{2\sqrt{{K}}}>\frac{d_h +\mu_h+\alpha N^*}{\rho}\Leftrightarrow \frac{r{\xi_h}}{2\sqrt{{K}}}-\frac{\alpha\sqrt{\left(\frac{r}{d}\right)^2-4\frac{K}{{\xi_h^2}}}}{2\rho}>\frac{d_h +\mu_h+\alpha \frac{r}{2d}}{\rho}.$$ Thus, Theorem \ref{th1:basics} provides a critical function of the hives population $\xi_h$ such that honeybee population can persist. In addition, notice that $\frac{r{\xi_h}}{2\sqrt{{K}}}-\frac{\alpha\sqrt{\left(\frac{r}{d}\right)^2-4\frac{K}{{\xi_h^2}}}}{2\rho}$ is an {increasing} function of $\xi_h$,  this implies that  the higher hives to brood ratio $\xi_h$, the better honeybee growth will be, and more likely persistent for  honeybees. This is supported by the empirical study of \cite{Schmickl:2007aa}.} \modifyc{In the following theorem, we provide theoretical results on the sufficient conditions that lead to the persistence and the  extinction of disease population or Varorra mites population.}\\

\begin{theorem}[Persistence and extinction of disease or mites]\label{th2:per-ext}
The following statements hold
\begin{itemize}
 \item If $N^*<\frac{d_m}{\alpha }$, then the total population of mite $N_m$ goes extinct, i.e.,
$$\limsup_{t\rightarrow\infty} N_m(t)=0$$where 
  system \eqref{newHBD} is attracted to  the mite-free invariant  set $MF=\{(S_h, I_h, S_m, I_m)\in\mathbb R^4_+:  S_m+ I_m=0\}$, and its dynamics is equivalent to the following two-D model \eqref{HV}
 \bae\label{HV}
\begin{array}{lcl}
S_h'&=& \frac{r (S_h+\rho I_h)^2}{\hat{K}+(S_h+\rho I_h)^2}-d_h S_h-S_h\frac{\beta_{h}I_h}{S_h+I_h}\\\\

I_h'&=&S_h\frac{\beta_{h}I_h}{S_h+I_h}-(d_h+\mu_h)I_h

\end{array}.
\eae

\item If $\frac{r}{2\sqrt{\hat{K}}}>d_h,\, \bar{N}_h^*<\frac{d_m}{\alpha c},\, \mbox{ and }S_h(0)>{S}_h^c$, then the total population of honeybees persists while the healthy mite population $S_m$ goes extinct, i.e.,
$$\limsup_{t\rightarrow\infty}S_m(t)=0$$ where the system \eqref{newHBD} is attracted to
 the healthy-mite-free invariant set $HMF=\{(S_h, I_h, S_m, I_m)\in\mathbb R^4_+:  S_m=0\}$ and its dynamics is equivalent to the following three-D system \eqref{HVMI}:
\bae\label{HVMI}
\begin{array}{lcl}
S_h'&=&\frac{r\left(S_h+\rho I_h\right)^2}{\hat{K}+\left(S_h+\rho I_h\right)^2}-d_h S_h-\frac{\beta_{h} S_hI_h}{S_h+I_h}-\modifyr{\frac{\hat{\beta}_{mh} S_hI_m}{S_h+I_h}}-\tilde{\beta}_{mh}S_hI_m-\alpha S_hI_m\\\\
I_h'&=&S_h\left[\frac{\beta_{h}I_h}{S_h+I_h}+\modifyr{\frac{\hat{\beta}_{mh} I_m}{S_h+I_h}}+\tilde{\beta}_{mh}I_m\right]-\alpha I_h I_m-(d_h+\mu_h)I_h\\\\
I_m'&=&c\alpha I_m\left[I_h +S_h-\frac{d_m +\mu_m}{c\alpha}\right].
\end{array}
\eae 

 \item Assume that $\frac{r}{2\sqrt{\hat{K}}}>d_h+\beta_h+\modifyr{\frac{\hat{\beta}_{mh} N^*}{\underline{N}^*_h}}+(\tilde{\beta}_{mh}+\alpha)(N^*-\underline{N}^*_h)\mbox{ and }S_h(0)>{S}_h^c$. Then the disease $I=cI_h+I_m$ persists if the inequality $\modifyr{\frac{\min\Big\{\beta_{h},\modifyr{c\hat{\beta}_{mh} }+c\tilde{\beta}_{mh}+c\alpha {S}^*_h\Big\}}{\max\Big\{(d_h+\mu_h),(d_m +\mu_m)\Big\}}\geq1}$ holds.
\item 
Assume that $\frac{r}{2\sqrt{\hat{K}}}>\frac{d_h+\mu_h+\alpha N^*}{\rho}\mbox{ and }N_h(0)>\underline{N}^c_h$. Then the disease $I=cI_h+I_m$ goes extinct if the following inequality holds
$\frac{\max\Big\{\beta_{h}+\frac{\hat{\beta}_{hm}N^*}{\underline{N}^*_h},\modifyr{c\hat{\beta}_{mh}}+c\tilde{\beta}_{mh}+c\alpha \bar{N}^*_h\Big\}}{\min\Big\{(d_h+\mu_h),(d_m +\mu_m)\Big\}}<1.$ Under this condition, the system \eqref{newHBD} is attracted to the virus-free invariant set $DF=\{(S_h, I_h, S_m, I_m)\in\mathbb R^4_+:  I_h+ I_m=0\}$ and its dynamics is equivalent to the following two-D model \eqref{HM}
 \bae\label{HM}
\begin{array}{lcl}
S_h'&=& \frac{rS_h^2}{\hat{K}+S_h^2}-d_h S_h-\alpha S_h S_m\\\\
S_m'&=&c\alpha S_hS_m-d_m S_m
\end{array}.
\eae

\end{itemize}

\end{theorem}

\noindent\textbf{Notes:} The results of the reduced dynamics in Theorem \ref{th2:per-ext} can be easily obtained by the theory of asymptotically autonomous systems \cite{Castillo-Chavez:1995aa}. The detailed proof of our results are provided in the last section. \\
According to Theorem \ref{th2:per-ext}, the condition $N^*=c\frac{\frac{r}{d_h}+\sqrt{\left(\frac{r}{d_h}\right)^2-4\hat{K}}}{2}<\frac{d_m}{\alpha}$ can lead to the extinction of the mite population.  Therefore, we can conclude that large values of the death rate of mites, $d_m$, can lead to the extinction of the whole colony; and large values of the death rate of mites $d_m$, small values of mite attacking rate, $\alpha$, and its energy conversion rate $c$, can lead to either the extinction of the whole mite population $N_m$ or the extinction of the healthy mite population $S_m$. Here we would like to point out that it is possible to have the persistence of \modifyb{virus infected} mites while the healthy mite goes extinct (see the resulting dynamics \eqref{HVMI} when the healthy mite goes extinct). In addition, the results of Theorem \ref{th2:per-ext} also suggest that: \textbf{1.} the persistence of the virus requires a large value for the disease transmission rate between adult honeybees, $\beta_h$, or that the disease transmission rate between honeybee brood and reproductive mites, $\tilde{\beta}_{mh}$; or small values of total death rates of honeybees, $d_h+\mu_h$, and mites, $d_m+\mu_m$; \textbf{2.} the extinction of the virus requires small values of all disease transmission rates, i.e., small values of $\beta_h, \hat{\beta}_{mh}, \tilde{\beta}_{mh}, \hat{\beta}_{hm}$; or large values of total death rates of honeybees and mites.\\


\noindent \modifyc{Theorem \ref{th2:per-ext} provides sufficient conditions that the full system \eqref{newHBD} reduces to the virus-free subsystem \eqref{HM},  the mite-free subsystem \eqref{HV}, and the healthy-mite-free subsystem \eqref{HVMI}. In the following three subsections, we explore the global dynamics of these subsystems. }\\
\subsection{Dynamics of the virus-free subsystem: only parasitism by mites}
\modifyb{Theorem \ref{th2:per-ext} in previous section suggests that either small values of all virus transmission rates or large values of total death rates of honeybees and mites can lead to the extinction of the virus infected honeybees and mites, which gives the following virus-free dynamics \eqref{HM}:}
$$\begin{array}{lcl}
S_h'&=& \frac{rS_h^2}{\hat{K}+S_h^2}-d_h S_h-\alpha S_h S_m\\\\
S_m'&=&c\alpha S_hS_m-d_m S_m
\end{array}$$
The dynamics of the virus-free dynamics \eqref{HM} (i.e., the dynamics of the parasitism interactions between honeybees and Varroa mites) can be summarized by the following theorem:\\

\begin{theorem}[Dynamics of the virus-free subsystem]\label{th5:DF} Let 
$H^*=\frac{d_m}{\alpha c},\,\, M^*=\frac{1}{\alpha}\left[\frac{rH^*}{\hat{K}+(H^*)^2}-d_h \right].$ \\
The virus-free subsystem \eqref{HM} can have one, three, or four equilibria. The existence and stability conditions for these equilibria are listed in Table \ref{t1:DF}.
\begin{table}[ht]
\centering
\begin{tabular}{|c|l|l|}
	\hline
 Equilibria &  Existence Condition        &   Stability Condition                \\
          \hline
$(0,0)$&         Always exists &Always locally stable \\\hline
$(\bar{N}^c_h,0)$& $\frac{r}{2\sqrt{\hat{K}}}>d_h$ &Saddle if $\bar{N}^c_h<\frac{d_m}{\alpha c}=H^*$; Source if $\bar{N}^c_h>\frac{d_m}{\alpha c}=H^*$\\\hline
 $(\bar{N}^*_h,0)$     &  $\frac{r}{2\sqrt{\hat{K}}}>d_h$ &Sink if $\bar{N}^*_h<\frac{d_m}{\alpha c}=H^*$; Saddle if $\bar{N}^*_h>\frac{d_m}{\alpha c}=H^*$\\\hline
$(H^*,M^*)$ & $\bar{N}^c_h<\frac{d_m}{\alpha c}=H^*<\bar{N}^*_h$ & Sink if $H^*>\sqrt{\hat{K}}$; Source if $H^*<\sqrt{\hat{K}}$. \\\hline
            \end{tabular}
\caption{The existence and stability of equilibrium for the virus-free subsystem \eqref{HM}, where 
$\bar{N}^c_h=\frac{\frac{r}{d_h}-\sqrt{\left(\frac{r}{d_h}\right)^2-4\hat{K}}}{2},\,\,\bar{N}^*_h=\frac{\frac{r}{d_h}+\sqrt{\left(\frac{r}{d_h}\right)^2-4\hat{K}}}{2}$ and 
$H^*=\frac{d_m}{\alpha c},\,\, M^*=\frac{1}{\alpha}\left[\frac{rH^*}{\hat{K}+(H^*)^2}-d_h \right].$}
\label{t1:DF}
\end{table}
The global dynamics of the virus-free subsystem \eqref{HM} can be summarized as follows:
\begin{enumerate}
\item The system \eqref{HM} converges to extinction $(0,0)$ for almost all initial conditions if $\frac{r}{2\sqrt{\hat{K}}}<d_h$ or $\frac{d_m}{\alpha c}<\bar{N}^c_h$.
\item If $\bar{N}^*_h<\frac{d_m}{\alpha c}$, depending on initial condition, the trajectory of \eqref{HM} converges to either $(0,0)$ or $(\bar{N}_h^*,0)$.
\item If $\bar{N}^c_h<\frac{d_m}{\alpha c}<\bar{N}^*_h$ 
 then Model  \eqref{HM} has a unique interior equilibrium $(H^*, M^*)$ which is locally asymptotically stable when $\frac{d_m}{\alpha c}>\sqrt{\hat{K}}$ and is a source when  $\frac{d_m}{\alpha c}<\sqrt{\hat{K}}$.
\end{enumerate}
\end{theorem}

\noindent\textbf{Notes:} 
Theorem \ref{th5:DF} provides us a global picture on the dynamics of the virus-free subsystem  \eqref{HM}, i.e., the honeybee colony only \modifyb{virus infected} with mites but not the virus.  By applying the results in \cite{Thieme:2009aa,Wang:2011aa}, we can conclude that the virus-free subsystem \eqref{HM} undergoes a subcritical Hopf-bifurcation at $\frac{d_m}{\alpha c}=\sqrt{\hat{K}}$. The subsystem  \eqref{HM}  has a unique \textbf{unstable} limit cycle around $(H^*,M^*)$ whenever $\frac{d_m}{\alpha c}<\sqrt{\hat{K}}$. In this case, the periodic orbits expand until it touches the stable manifold of the boundary equilibrium $(\bar{N}^c_h,0)$ which leads to the extinction of both honeybees and the parasitic mites. We refer to this phenomena as a \emph{catastrophic event} which could be linked to CCD. Our theoretical results also suggest that a small death rate for mites and a large parasitism rate can destabilize the system.\\

\begin{figure}
\begin{center}
\includegraphics[width=150mm]{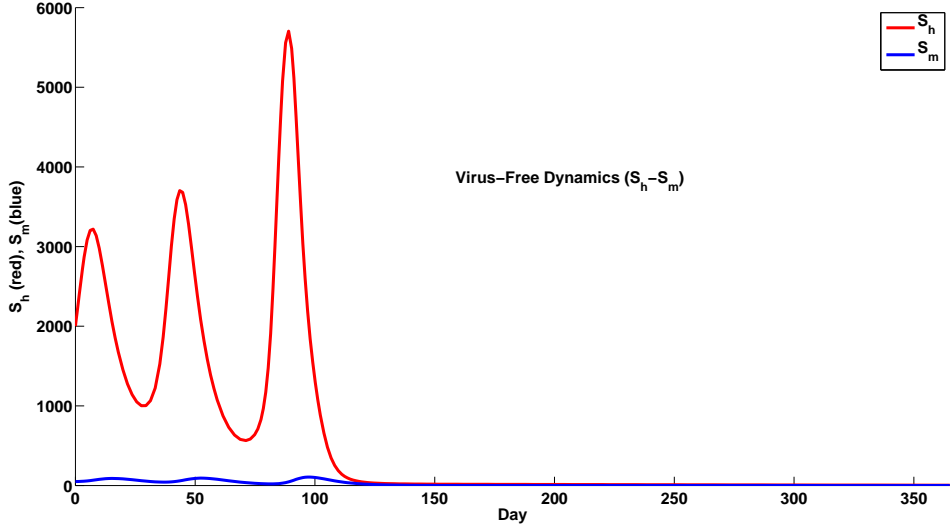}
\end{center}
\caption{Time series \modifyr{(in days)} of Model \eqref{HM} when $r=1500,\,\alpha = 0.005,\,c = 0.01,\, d_h=0.01,\, d_m=0.1$: population of honeybees is in red while Varroa mites is in blue.}
\label{fig1:DF}
\end{figure}

\noindent\textbf{Linking to CCD:}  To illustrate the  \emph{catastrophic event}, we use reasonable parameters from \cite{Sumpter:2004aa,Ratti:2012aa}.
Let the reproduction of egg per day during summer be $r=1500$; and the population size of the honeybee colony at which the birth rate is
half of the maximum possible rate  be $\sqrt{\hat{K}}=2000$; the natural death rate of honeybees is $d_h= 0.01$; the parasitism rate is $\alpha = 0.005$; the energy conversion rate is $c = 0.01$; and the natural death rate of mites is $d_m = 0.1$. This set of parameter values gives $\frac{d_m}{\alpha c}<\sqrt{\hat{K}}$ which implies that a \emph{catastrophic event} will occur (see Figure \ref{fig1:DF}; the population of honeybees is in red and collapses around \modifyr{time=200 day}).\\

\modifyb{\noindent\textbf{Stochastic effects and oscillations:} Theorem \ref{th5:DF} implies that if the inequality $\frac{d_m}{\alpha c}<\sqrt{\hat{K}}$ holds, then the virus-free subsystem  \eqref{HM} has a unique \textbf{unstable} limit cycle around $(H^*,M^*)$ where for all initial conditions either the system converges to $(0,0)$ quickly or the system experiences the expanding oscillations leading to the extinction eventually. The oscillating extinction in the later type is driven by the deterministic dynamics. The extinction fate of the system cannot be prevented by introducing stochastic effects, however, introduced stochastic effects may cause the system going to extinction more quickly without expanding oscillations.  }\\

\modifyr{Note that $\hat{K}=\frac{K}{\xi_h^2}$ and $\alpha=\hat{\alpha} (1-\xi_h)(1-\xi_m)$ where $\xi_h,\xi_m$ are ratio of adult bees and ratio of phoretic stage of Varroa mites in honeybee colony, respectively. The \emph{catastrophic event} occurs when 
$$\frac{d_m}{\alpha c}<\sqrt{\hat{K}} \Leftrightarrow \frac{d_m}{\hat{\alpha}(1-\xi_m) c\sqrt{K}}<\frac{ (1-\xi_h)}{\xi_h}\Leftrightarrow \frac{\xi_h}{ (1-\xi_h)}<\frac{\hat{\alpha}(1-\xi_m) c\sqrt{K}}{d_m}.$$ This inequality provides a critical low hive to brood ratio that can destabilize the system and cause the sudden extinction of honeybees.}\\

\subsection{Dynamics of the mite-free subsystem: only virus infections}
\modifyb{According to Theorem \ref{th2:per-ext}, if the honeybee population is too small, e.g., $N^*<\frac{d_m}{\alpha c}$, then the dynamics of \eqref{newHBD} is  equivalent to the following mite-free dynamics \eqref{HV} 
$$\begin{array}{lcl}
S_h'&=& \frac{r (S_h+\rho I_h)^2}{\hat{K}+(S_h+\rho I_h)^2}-d_h S_h-S_h\frac{\beta_{h}I_h}{S_h+I_h}\\\\
I_h'&=&S_h\frac{\beta_{h}I_h}{S_h+I_h}-(d_h+\mu_h)I_h
\end{array}.$$
To continue studying the dynamics of the mite-free dynamics \eqref{HV} , we define $a=\frac{1}{\frac{\beta_h}{d_h+\mu_h}-1}$ as the ratio of the susceptible honeybee population to the virus infected honeybee population; $\mathcal R_0^V=\frac{\beta_h}{d_h+\mu_h}$ as the basic reproduction number, i.e., the number of secondary cases which one case would produce in a completely susceptible population;  $\tilde{d}=(a+1) d_h+\mu_h=d_h\left(\frac{\mathcal R_0^V}{\mathcal R_0^V-1}\right)+\mu_h$ as the updated average death of the honeybee due to virus infections. In addition, we let $S_h^k= a I_h^k, k=1,2$ and
$$\begin{array}{lll}
I_h^{1}=\frac{\frac{r}{\tilde{d}}-\sqrt{\left(\frac{r}{\tilde{d}}\right)^2-4\frac{\hat{K}}{(a+\rho)^2}}}{2},&&\,\,I_h^{2}=\frac{\frac{r}{\tilde{d}}+\sqrt{\left(\frac{r}{\tilde{d}}\right)^2-4\frac{\hat{K}}{(a+\rho)^2}}}{2}.
\end{array}$$
The dynamics of the  mite-free system \eqref{HV} can be summarized by the following theorem:}

\begin{theorem}[Dynamics of the mite-free subsystem]\label{th3:MF}
The mite-free subsystem \eqref{HV} can have one, three, or five equilibria. 
The existence and stability conditions for these equilibria are listed in Table \ref{t2:MF}.
\begin{table}[ht]
\centering
\begin{tabular}{|c|l|l|}
	\hline
 Equilibria &  Existence Condition for Existence       &   Stability Condition                \\
          \hline
$(0,0)$&         Always exists &Always locally stable \\\hline
$(\bar{N}^c_h,0)$& $\frac{r}{2\sqrt{\hat{K}}}>d_h$ &Saddle if $\mathcal R_0^V<1$; Source if $\mathcal R_0^V>1$\\\hline
 $(\bar{N}^*_h,0)$     &  $\frac{r}{2\sqrt{\hat{K}}}>d_h$ &Sink if $\mathcal R_0^V<1$; Saddle if $\mathcal R_0^V>1$\\\hline
$(S_h^1,I_h^1)$ & $\mathcal R_0^V>1$ and  $\frac{r}{2\sqrt{\hat{K}}}>\frac{\tilde{d}}{a+\rho}$& Always a saddle. \\\hline
  $(S_h^2,I_h^2)$ &  $\mathcal R_0^V>1$ and  $\frac{r}{2\sqrt{\hat{K}}}>\frac{\tilde{d}}{a+\rho}$  & Always locally asymptotically stable   \\\hline
            \end{tabular}
\caption{The existence and stability of equilibrium for the mite-free subsystem \eqref{HV}. We have $\mathcal R_0^V=\frac{\beta_h}{d_h+\mu_h}$,\, $a=\frac{1}{\mathcal R_0^V-1},\,\,\tilde{d}=(a+1) d_h+\mu_h=d_h\left(\frac{\mathcal R_0^V}{\mathcal R_0^V-1}\right)+\mu_h,$\, $\bar{N}^c_h=\frac{\frac{r}{d_h}-\sqrt{\left(\frac{r}{d_h}\right)^2-4\hat{K}}}{2},\,\,\bar{N}^*_h=\frac{\frac{r}{d_h}+\sqrt{\left(\frac{r}{d_h}\right)^2-4\hat{K}}}{2},$ and $I_h^{1}=\frac{\frac{r}{\tilde{d}}-\sqrt{\left(\frac{r}{\tilde{d}}\right)^2-4\frac{\hat{K}}{(a+\rho)^2}}}{2},\,\,I_h^{2}=\frac{\frac{r}{\tilde{d}}+\sqrt{\left(\frac{r}{\tilde{d}}\right)^2-4\frac{\hat{K}}{(a+\rho)^2}}}{2},\,\, S_h^k= a I_h^k, k=1,2.$}
\label{t2:MF}
\end{table}
In addition, the global dynamics of the mite-free subsystem \eqref{HV} can be summarized as follows:
\begin{enumerate}
\item The trajectory of \eqref{HV} converges to extinction $(0,0)$ for all initial conditions in $\mathbb R^2_+$ if one of the following conditions hold:
\begin{itemize}
\item $\frac{r}{2\sqrt{\hat{K}}}<d_h$, or 
\item $\mathcal R_0^V>1 \mbox{ and } \,d_h< \frac{r}{2\sqrt{\hat{K}}}<\frac{\tilde{d}}{a+\rho}$.
\end{itemize}

\item The trajectory of \eqref{HV} converges to either $(0,0)$ or $(\bar{N}_h^*,0)$ for almost all initial conditions in $\mathbb R^2_+$ if the inequalities $\mathcal R_0^V<1$ and $\frac{r}{2\sqrt{\hat{K}}}>d_h$.
\item The trajectory of \eqref{HV} converges to either $(0,0)$ or $(S^2_h,I^2_h)$ for almost  all initial conditions in $\mathbb R^2_+$ if   the inequalities $\mathcal R_0^V>1,\,\,\frac{r}{2\sqrt{\hat{K}}}>\max\Big\{d_h, \frac{\tilde{d}}{a+\rho}\Big\}$ hold.

\end{enumerate}
\end{theorem}

\noindent\textbf{Notes:} Theorem \ref{th3:MF} implies that the mite-free subsystem \eqref{HV}  has relatively simple dynamics, i.e., no limit cycle. The results show the following interesting findings:
\begin{enumerate}
\item Honeybees can persist with proper initial conditions if the virus transmission rate among honeybees $\beta_h$ is not large, i.e., $\mathcal R_0^V<1$.
\item Both honeybees and the virus can coexist if $\beta_h$ is in the medium range, i.e. $\mathcal R_0^V>1$ and $\frac{\tilde{d}}{a+\rho}<\frac{r}{2\sqrt{\hat{K}}}$
\item However, the large virus transmission rate among honeybees $\beta_h$ can drive honeybees to extinction. This occurs when the inequalities $\mathcal R_0^V>1 \mbox{ and } \,d_h< \frac{r}{2\sqrt{\hat{K}}}<\frac{\tilde{d}}{a+\rho}$ hold.\\

\end{enumerate}


\subsection{Dynamics of the healthy-mite-free subsystem}
 According to Theorem \ref{th2:per-ext}, if $\frac{r}{2\sqrt{\hat{K}}}>d_h,\, \bar{N}_h^*<\frac{d_m}{\alpha c},\, \mbox{ and }S_h(0)>{S}_h^c$, then the total population of honeybees persists while the healthy mite population $S_m$ goes extinct, i.e.,
$$\limsup_{t\rightarrow\infty}S_m(t)=0$$ where the system \eqref{newHBD} is attracted to
 the healthy-mite-free invariant set $HMF=\{(S_h, I_h, S_m, I_m)\in\mathbb R^4_+:  S_m=0\}$ and its dynamics is equivalent to the following three-D system \eqref{HVMI}:
$$\begin{array}{lcl}
S_h'&=&\frac{r\left(S_h+\rho I_h\right)^2}{\hat{K}+\left(S_h+\rho I_h\right)^2}-d_h S_h-\frac{\beta_{h} S_hI_h}{S_h+I_h}-\modifyr{\frac{\hat{\beta}_{mh} S_hI_m}{S_h+I_h}}-\tilde{\beta}_{mh}S_hI_m-\alpha S_hI_m\\\\
I_h'&=&S_h\left[\frac{\beta_{h}I_h}{S_h+I_h}+\modifyr{\frac{\hat{\beta}_{mh} I_m}{S_h+I_h}}+\tilde{\beta}_{mh}I_m\right]-\alpha I_h I_m-(d_h+\mu_h)I_h\\\\
I_m'&=&c\alpha I_m\left[I_h +S_h-\frac{d_m +\mu_m}{c\alpha}\right]
\end{array}.$$
 Let
$$\begin{array}{lcl}
f_1(I_h)&=&\frac{\frac{r\left(\frac{d_m+\mu_m}{c\alpha}-I_h+\rho I_h\right)^2}{\hat{K}+\left(\frac{d_m+\mu_m}{c\alpha}-I_h+\rho I_h\right)^2}-d_h \frac{d_m+\mu_m}{c\alpha}-\mu_h I_h}{\alpha(\frac{d_m+\mu_m}{c\alpha})}\\\\
 f_2(I_h)&=&\frac{\left(\frac{d_m+\mu_m}{c\alpha}-I_h\right)\modifyr{\left[\frac{\beta_{h}}{d_m +\mu_m}-(d_h+\mu_h)\right]I_h}}{\left[\alpha+\modifyr{\frac{c\alpha\hat{\beta}_{mh}}{d_m+\mu_m}}+\tilde{\beta}_{mh}\right] I_h-\frac{d_m+\mu_m}{c\alpha}\left(\modifyr{\frac{c\alpha\hat{\beta}_{mh}}{d_m+\mu_m}}+\tilde{\beta}_{mh}\right)}.\end{array}$$
The dynamics of the healthy-mite-free subsystem \eqref{HVMI} can be summarized by the following theorem:\\

\modifyr{\begin{theorem}[Dynamics of the healthy-mite-free subsystem]\label{th4:HVMI} 
If $\frac{r}{2\sqrt{\hat{K}}}>d_h$ and $\bar{N}_h^*<\frac{d_m+\mu_m}{c\alpha}$, then the population of \modifyb{virus infected} mites goes extinct in the subsystem \eqref{HVMI}, i.e., $$\limsup_{t\rightarrow\infty}I_m(t)= 0$$  which reduces to the following mite-free model \eqref{HV}. In addition, the following statements hold
\begin{enumerate}
\item If $(S_h, I_h,I_m)$ is an interior equilibrium of the healthy-mite-free subsystem \eqref{HVMI}, then $I_h$ is a positive intercept of $f_1(I_h)$ and $f_2(I_h)$ subject to $0<I_h<\frac{d_m+\mu_m}{c\alpha}$, $S_h=\frac{d_m+\mu_m}{c\alpha}-I_h$ and $I_m=f_1(I_h)$.
\item The healthy-mite-free subsystem \eqref{HVMI} has no interior equilibrium if the inequality $ rc\alpha<d_h (d_m+\mu_m)$ holds.
\item Assume that $\frac{r}{2\sqrt{\hat{K}}}>\frac{d_h+\mu_h+\alpha N^*}{\rho}\mbox{ and }N_h(0)>\underline{N}^c_h$. Then both virus infected honeybee population $I_h$ and virus infected mites $I_m$ persist if the inequalities $\bar{N}_h^*>\frac{d_m+\mu_m}{c\alpha}$ and $\mathcal R_0^V=\frac{\beta_h}{d_h+\mu_h}<1$ hold.
\end{enumerate}
\end{theorem}}

\noindent\textbf{Notes:}  Theorem \ref{th4:HVMI} implies that, under the condition of $\bar{N}_h^*<\frac{d_m+\mu_m}{c\alpha}$, the \modifyb{virus infected} mite population $I_m$ goes extinct in the healthy-mite-free subsystem \eqref{HVMI} which reduces to the mite-free subsystem \eqref{HV} that we studied in the previous subsection. In addition, Theorem \ref{th4:HVMI} shows that the subsystem \eqref{HVMI} has no interior equilibrium if the inequality $rc\alpha<d_h (d_m+\mu_m)$ holds. Therefore, we could expect the extinction of $I_m$ for small values of $r, c, \alpha$ and large values of $d_h, d_m, \mu_m$. This has been confirmed by numerical simulations. The population of honeybees and virus infected mites in \eqref{HVMI} experiences sudden collapse when we (1) increase the values of $c,\alpha, \hat{K},$ and the related virus transmission rates, or (2) decrease the values of $d_h, \rho, r, \mu_m$. The biological implication for this dynamics is that increasing or decreasing the values of these parameters destabilizes the system and generates fluctuated dynamics. The destabilizing effects generate unstable oscillations. The amplitudes of oscillations increase until they touch the stable manifold of the extinction equilibrium, which cause the collapse of the whole colony. \modifyr{The destabilizing effects of $c,\alpha, \hat{K}, d_m$ can be explained through the dynamics of the virus free subsystem \eqref{HM} that we have studied in  Theorem \ref{th5:DF}. 
\begin{itemize}
\item Decreasing the values of $c,\alpha, \hat{K}$ can stabilize the system; small values of  $c,\alpha$ can cause the extinction of the \modifyb{virus infected} mite population $I_m$, and lead to the coexistence of $S_h$ and $I_h$.
\item Increasing the value of $\mu_h$ can stabilize the system but large values of $\mu_h$ can cause extinction of the whole colony due to the initial oscillations. 
\item Decreasing $\mu_m$ can destabilize the system; while increasing it can stabilize the system; large values of $\mu_m$ can lead to the extinction of $I_m$ and the persistence  of $S_h, I_h$.
\item Decreasing the value of $\rho$ could destabilize the system, thus causing the extinction of the colony.
\item Increasing the virus transmission rates (i.e., $\beta_h, \hat{\beta}_{hm}, \tilde{\beta}_{mh},\hat{\beta}_{hm}$) can stabilize the system, while decreasing their values can destablize the system and cause the extinction of all species.\\
\end{itemize}}

\begin{figure}[ht]
\centering
   \includegraphics[scale =.4] {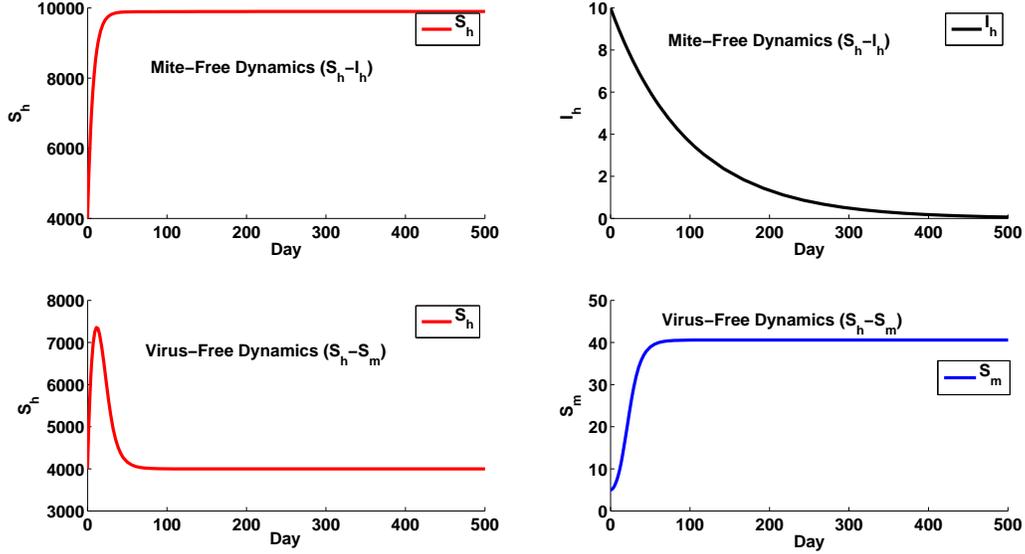}
  \caption{Population dynamics of the subsystems of the the honeybee-mite-virus model \eqref{newHBD} when $r=1500; K=1000000; \rho=0.9; d_h=.15; \mu_h=0.1; \alpha=0.005; d_m=0.1; \mu_m=0.01;  c=0.005; \beta_h=.24; \hat{\beta}_{mh} = 0.03; \tilde{\beta}_{mh} = .005; \hat{\beta}_{hm} = 0.03$. The left figure in the first row is the healthy honeybee population (the red curve) and the right figure in the first row is the virus infected honeybee population (the black curve) in the mite-free subsystem \eqref{HM} when  $S_h(0)=4001, I_h(0)=10$.  The left figure in the second row is the healthy honeybee population (the red curve) and the right figure in the second row is the healthy mite population (the blue curve) in the virus-free subsystem \eqref{HV} when$S_h(0)=4001,S_m(0)=5$. 
  }\label{fig2:comparison1}
\end{figure}

\begin{figure}[ht]
\centering
   \includegraphics[scale =.4] {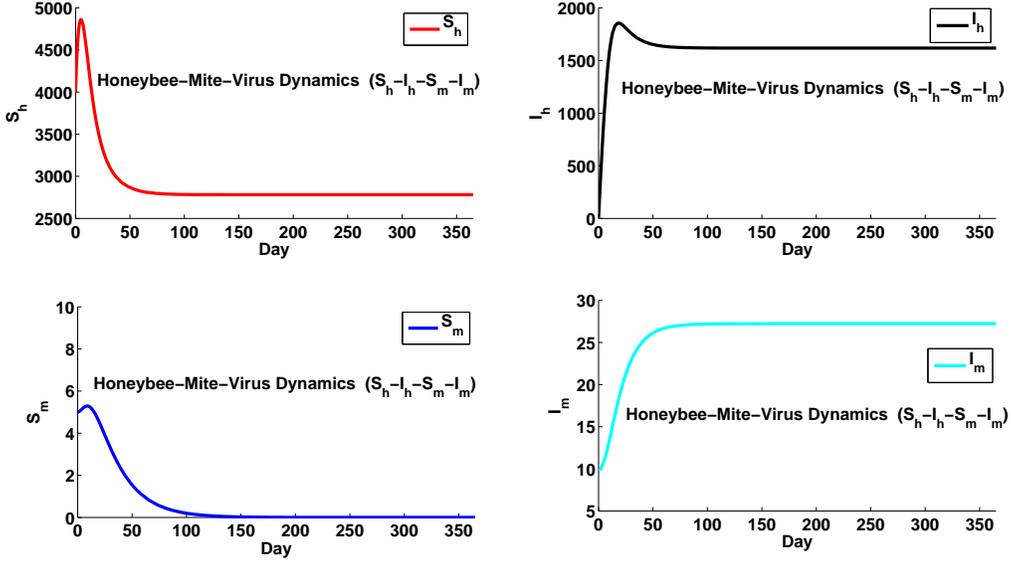}
  \caption{Population dynamics of the honeybee-mite-virus model \eqref{newHBD} when $r=1500; K=1000000; \rho=0.9; d_h=.15; \mu_h=0.1; \alpha=0.005; d_m=0.1; \mu_m=0.01;  c=0.005; \beta_h=.24; \hat{\beta}_{mh} = 0.03; \tilde{\beta}_{mh} = .005; \hat{\beta}_{hm} = 0.03$ and $S_h(0)=4001, I_h(0)=10, S_m(0)=5, I_m(0)=10.$
The healthy honeybee population $S_h$ is in red; the \modifyb{virus infected} honeybee population $I_h$ is in black; the healthy mite population $S_m$ is in blue; and the \modifyb{virus infected} mite population $I_m$ is in cyan.}\label{fig2:comparison2}
\end{figure}

\noindent\textbf{Synergistic effects of parasitic mites and virus infections:}  \modifyr{If there is no mites in the system, according to Theorem \ref{th3:MF}, the mite-free subsystem \eqref{HV} reduces to the only healthy honeybee population, i.e., the disease goes extinct whenever the initial population of honeybees is above $\bar{N}^c_h$, $\frac{r}{2\sqrt{\hat{K}}}>\frac{d_h+\mu_h+\alpha N^*}{\rho}$, and the basic reproduction number $\mathcal R_0^V=\frac{\beta_h}{d_h+\mu_h}<1$ (see two figures in the first row of Figure \ref{fig2:comparison1} where virus infected honeybees go extinct (the black curve in right) and the healthy honeybees persist (the red curve in red)). However, when the virus infected mites are in the colony, i.e., the healthy-mite-free subsystem \eqref{HVMI}, both disease and mites can persist under proper conditions (see Figure \ref{fig2:comparison2} where virus infected honeybees (the black curve), the healthy honeybees (the red curve), and the virus infected mites persist (the cyan curve) with the healthy mites going extinct (the blue curve). For example, the synergistic effects of parasitic mites and virus infections has been illustrated in Figure \ref{fig2:comparison1}-\ref{fig2:comparison2} when $r=1500; K=1000000; \rho=0.9; d_h=.15; \mu_h=0.1; \alpha=0.005; d_m=0.1; \mu_m=0.01;  c=0.005; \beta_h=.24; \hat{\beta}_{mh} = 0.03; \tilde{\beta}_{mh} = .005; \hat{\beta}_{hm} = 0.03$. These are reasonable parameter values derived from \cite{Sumpter:2004aa,Ratti:2012aa}. }\\

\subsection{Dynamics of the full system}
Recall that the full system \eqref{newHBD} of honeybee-mite-virus interactions can be described by the following set of equations:
$$\begin{array}{lcl}
S_h'&=& \frac{r\left(S_h+\rho I_h\right)^2}{\hat{K}+\left(S_h+\rho I_h\right)^2}-d_h S_h-\frac{\beta_{h} S_hI_h}{S_h+I_h}-\modifyr{\frac{\hat{\beta}_{mh} S_hI_m}{S_h+I_h}}-\tilde{\beta}_{mh}S_hI_m-\alpha S_h (S_m+I_m)\\\\
I_h'&=&S_h\left[\frac{\beta_{h}I_h}{S_h+I_h}+\modifyr{\frac{\hat{\beta}_{mh}I_m}{S_h+I_h}}+\tilde{\beta}_{mh}I_m\right]-\alpha I_h (S_m+I_m)-(d_h+\mu_h)I_h\\\\
S_m'&=&S_m\left[c\alpha S_h-\frac{\hat{\beta}_{hm}I_h}{S_h+I_h}-d_m\right]\\\\
I_m'&=&c\alpha\left[I_h (S_m+I_m)+S_hI_m\right]+\frac{\hat{\beta}_{hm}I_hS_m}{S_h+I_h}-(d_m +\mu_m)I_m
\end{array}.
$$
The results from the previous section provide us a complete picture of the dynamics of the subsystems of the full system \eqref{newHBD}. In this subsection, we explore the dynamics of the full system as the following theorem. \\
\modifyr{\begin{theorem}[The persistence of honeybees]\label{th6:global} Assume that  $\frac{r}{2\sqrt{\hat{K}}}>d_h$. If $N^*<\frac{d_m}{\alpha }$ and $\mathcal R_0^V=\frac{\beta_h}{d_h+\mu_h}<1$, the full system \eqref{newHBD} converges to
the disease-mite-free set $DMF=\{(S_h, I_h, S_m, I_m)\in\mathbb R^4_+:  S_m+I_h+ I_m=0\}$ where the system \eqref{newHBD} is reduced to the following one-D system \eqref{H}:
\bae\label{H}
\begin{array}{lcl}
S_h'&=& \frac{rS_h^2}{\hat{K}+S_h^2}-d_h S_h
\end{array}
\eae whose dynamics can be summarized as follows:
\begin{enumerate}
\item If the inequality $\frac{r}{2\sqrt{\hat{K}}}<d_h$ holds, then \eqref{H} converges to 0.
\item If  the inequalities $\frac{r}{2\sqrt{\hat{K}}}>d_h$ and $S_h(0)> \bar{N}^c_h$ hold, then  \eqref{H}  converges to $\bar{N}^*_h$. If the initial condition falls below $\bar{N}^c_h$, i.e., $S_h(0)< \bar{N}^c_h$, then  \eqref{H} also converges to 0.
\end{enumerate}
Moreover, the following statements hold
\begin{enumerate}
\item If the inequalities $\frac{r}{2\sqrt{\hat{K}}}>\max\Big\{d_h+\beta_h+\modifyr{\frac{\hat{\beta}_{mh} N^*}{\underline{N}^*_h}}+(\tilde{\beta}_{mh}+\alpha)(N^*-c\underline{N}^*_h),\, \frac{d_h +\mu_h+\alpha N^*}{\rho}\Big\}$ \mbox{ with } $ N_h(0)\geq S_h(0)>S_h^c $ hold, then $S_h$ is persistent, i.e., $\liminf_{t\rightarrow\infty} S_h(t)\geq S^*_h.$
\item Assume that $\frac{r}{2\sqrt{\hat{K}}}>d_h+\beta_h+\modifyr{\frac{\hat{\beta}_{mh} N^*}{\underline{N}^*_h}}+(\tilde{\beta}_{mh}+\alpha)(N^*-\underline{N}^*_h)\mbox{ with }S_h(0)>{S}_h^c$ and $\bar{N}^*_h<\frac{d_m}{\alpha c}$ hold. Then the disease $I=cI_h+I_m$ persists if the inequality $\modifyr{\frac{\min\Big\{\beta_{h},\modifyr{c\hat{\beta}_{mh} }+c\tilde{\beta}_{mh}+c\alpha \bar{N}^*_h\Big\}}{\max\Big\{(d_h+\mu_h),(d_m +\mu_m)\Big\}}\geq1}$ holds.
\item The full system \eqref{newHBD} has no interior equilibrium if one of the following inequalities hold:\\\\
$\frac{\beta_{h}}{d_h+\mu_h}>\frac{\hat{\beta}_{hm}}{d_m+\hat{\beta}_{hm}},\,\frac{\hat{\beta}_{mh}}{\frac{\alpha\mu_m \hat{\beta}_{hm}}{(d_m+\hat{\beta}_{hm})}+c\alpha \hat{\beta}_{mh}-\tilde{\beta}_{mh} \hat{\beta}_{hm}}>\frac{d_m+\mu_m}{c\alpha(d_m+\hat{\beta}_{hm}+\mu_m)}.$ Or\\\\
$\frac{\alpha\mu_m \hat{\beta}_{hm}}{\hat{\beta}_{hm}+d_m}> \tilde{\beta}_{mh} \hat{\beta}_{hm}-c\alpha,\,0<\frac{\frac{\beta_{h}}{d_h+\mu_h}-\frac{\hat{\beta}_{hm}}{d_m+\hat{\beta}_{hm}}}{ \beta_{h}(d_m+\hat{\beta}_{hm})(d_h+\mu_h)}<d_m.$ Or\\\\
$\frac{\beta_{h}}{d_h+\mu_h}<\frac{\hat{\beta}_{hm}}{d_m+\hat{\beta}_{hm}},\,\,0<\frac{\hat{\beta}_{mh}(d_m+\hat{\beta}_{hm})}{\frac{\alpha\mu_m \hat{\beta}_{hm}}{(d_m+\hat{\beta}_{hm})}+c\alpha \hat{\beta}_{mh}-\tilde{\beta}_{mh} \hat{\beta}_{hm}}<\frac{d_m}{c\alpha}.$ Or\\\\
$\frac{\alpha\mu_m \hat{\beta}_{hm}}{\hat{\beta}_{hm}+d_m}< \tilde{\beta}_{mh} \hat{\beta}_{hm}-c\alpha,\,\,\frac{\frac{\beta_{h}}{d_h+\mu_h}-\frac{\hat{\beta}_{hm}}{d_m+\hat{\beta}_{hm}}}{ \beta_{h}(d_m+\hat{\beta}_{hm})(d_h+\mu_h)}>\frac{(d_m+\mu_m)(d_m+\hat{\beta}_{hm})}{d_m+\hat{\beta}_{hm}+\mu_m}.$
\item Assume that $\frac{r}{2\sqrt{\hat{K}}}>\frac{d_h+\mu_h+\alpha N^*}{\rho}\mbox{ and }N_h(0)>\underline{N}^c_h$. Then the total mite population $N_m=S_m+I_m$ persists if the inequalities $\bar{N}_h^*>\frac{d_m+\mu_m}{c\alpha}$ and $\mathcal R_0^V=\frac{\beta_h}{d_h+\mu_h}<1$ hold.\\
\end{enumerate}
\end{theorem}}

\noindent\textbf{Notes:} \modifyb{Theorem \ref{th6:global} along with Theorem \ref{th1:basics} - \ref{th5:DF}, we can conclude that the extinction of disease occurs when all values of all disease transmission rates, $\beta_h, \hat{\beta}_{mh}, \tilde{\beta}_{mh}, \hat{\beta}_{hm}$ are small; with the consequence that  the full system  \eqref{newHBD} converges to either $(0,0,0,0)$ or $(\bar{N}_h^*,0,0,0)$ when $\mathcal R_0^M=\frac{\bar{N}_h^*}{H^*}<1$ while  \eqref{newHBD} converges to  either $(0,0,0,0)$ or $(H^*, 0,M^*,0)$ when $1<\mathcal R_0^M<\frac{\bar{N}_h^*}{S^c_h}$. The persistence of disease or mites indicates the persistence of honeybees even though the population may be low under the influence of disease or mites. Theorem \ref{th6:global} provides a summary on sufficient conditions when honeybees can persist in the full \eqref{newHBD} alone, with mites, or with disease. The item 4 of Theorem \ref{th6:global}  is consistent with the results from Theorem \ref{th4:HVMI} regarding the synergistic effects of parasitic mites and virus infections: If there is no mites in the system, according to Theorem \ref{th3:MF}, the mite-free subsystem \eqref{HV} reduces to the only healthy honeybee population, however, in the presence of mites, both disease and mites in the honeybee-mite-virus system \eqref{newHBD} can persist under proper conditions (see Figure \ref{fig2:comparison1}-\ref{fig2:comparison2} for more details).}\\

The dynamics of the full system \eqref{newHBD} can be extremely complicated. We are unable to obtain an explicit form of the interior equilibrium and the related stability. We perform a series of numerical simulations to explore how different parameters affect the population dynamics. The effects of $r, c, \alpha$, $d_h, d_m, \mu_m$, $\mu_h, \rho$ and virus transmission rates are similar to our observations for the subsystem  \eqref{HVMI}. More specifically, we have the following observations:
\begin{enumerate}
\item Effects of $r$: Increasing $r$ can stabilize the system, but increasing it too large can drive healthy mites $S_m$ to extinction while the population of \modifyb{virus infected} mites $I_m$ increases. Decreasing the values of $r$ can destablize the system, and cause the extinction of mites. And too small values can cause the whole colony to become extinct.
\item Effects of $c$:  Increasing can destabilize the system and cause extinction of all species. Decreasing its value can stabilize the system but too small values can drive the extinction of mites. 
\item Effects of $\alpha$: Increasing can destabilize the system and cause extinction of all species while decreasing can stabilize the system and increase the healthy honeybee population. Small values of $\alpha$ can drive the mite population to extinction, and too small values can drive both mites and virus to extinction, and only healthy honeybees are left.

\item Effects of $d_m$:  Increasing $d_m$ can stabilize the system, and drive $S_m$ to extinction first. Increasing it further can lead to the extinction of mites, and the system approaches the limiting mite-free system \eqref{HV}.  Decreasing can destabilize the system and drive the extinction of the virus. Too small values can cause the whole colony's extinction.
\item Effects of $d_h$: Increasing can stabilize the system. Large values can drive the virus extinctions, however, extremely large values can lead to the extinction of the whole colony.
Decreasing can destabilize the system which may cause the extinction of the colony under certain conditions.

\item Effects of $\mu_m$: Increasing can drive virus extinction. Extremely large values can lead to the extinction of the whole colony.
Decreasing can destabilize the system. Small values can drive the healthy mites extinct, and extremely small values may cause the extinction of the whole colony.

\item Effects of $\mu_h$: Increasing can cause the extinction of virus. Decreasing can stabilize the system, and small values can drive healthy mites to extinction. 

\item	Effects of virus transmission rates: Increasing can destablize the system and cause the extinction of healthy mites $S_m$, while extremely large values may drive all populations extinct. Decreasing can drive the extinction of virus.\\
\end{enumerate}

\subsection{Mechanisms of collapse dynamics and synergistic effects}
\modifyc{Let $r=1500; K=1000000; \rho=0.9; d_h=.15; \mu_h=0.1; \alpha=0.005; d_m=0.1; \mu_m=0.01;  c=0.005; \beta_h=.24; \hat{\beta}_{mh} = 0.03; \tilde{\beta}_{mh} = .005; \hat{\beta}_{hm} = 0.03$ (Figure \ref{fig2:comparison1}-\ref{fig2:comparison2}) and $r = 1500; \rho=0.9; \hat{K}= 1600001; d_h =0.15; \mu_h=0.1; \alpha=0.05; c = 0.005; d_m=0.1;\mu_m=0.01; \beta_h=0.3; \hat{\beta}_{mh}=0.08; \tilde{\beta}_{mh}=0.001;\hat{\beta}_{hm}=0.03$ (Figure \ref{fig3:comparison1}-\ref{fig3:comparison2}). These are reasonable parameter values derived from \cite{Sumpter:2004aa,Ratti:2012aa}. We use these two sets of parameters as illustrations to explore the synergistic effects of parasite mites and virus infections as well as potential mechanisms linking to CCD (see Figure \ref{fig1:DF} and Figure \ref{fig3:comparison1}-\ref{fig3:comparison2}). These comparisons suggest the following:\\
\begin{enumerate}
\item Synergistic effects of parasitic mites and virus infections: Based on the two sets of parameters, we have the following two typical scenarios:
\begin{enumerate}
\item Under the parameter values of $r=1500; K=1000000; \rho=0.9; d_h=.15; \mu_h=0.1; \alpha=0.005; d_m=0.1; \mu_m=0.01;  c=0.005; \beta_h=.24; \hat{\beta}_{mh} = 0.03; \tilde{\beta}_{mh} = .005; \hat{\beta}_{hm} = 0.03$:\\\\
 If there is no mites, the mite-free system \eqref{HV} (i.e., the dynamics of the healthy honeybee and the virus infected honeybee) converges to only the healthy honeybee with virus infected honeybees go extinct (see the first row of Figure \ref{fig2:comparison1}). \\\\
 If there is no virus, the virus-free system \eqref{HM} (i.e., the parasitism dynamics of the healthy honeybee and the healthy mites) converges to a stable equilibrium where both the healthy honeybee and the healthy mites can persist (see the second row of Figure \ref{fig2:comparison1}).\\\\
 However, if honeybees, mites, and virus are all presented in the system (i.e., the full system \eqref{newHBD}), then both virus infected honeybees (black curve) and virus infected mites (cyan curve) can persist (see Figure \ref{fig2:comparison2}).\\\\
 This implies that the presence of mite population can promote the persistence of disease. \\
  \item Under the parameter values of $r = 1500; \rho=0.9; \hat{K}= 1600001; d_h =0.15; \mu_h=0.1; \alpha=0.05; c = 0.005; d_m=0.1;\mu_m=0.01; \beta_h=0.3; \hat{\beta}_{mh}=0.08; \tilde{\beta}_{mh}=0.001;\hat{\beta}_{hm}=0.03$:\\\\
 If there is no mites, the mite-free system \eqref{HV} (i.e., the dynamics of the healthy honeybee and the virus infected honeybee) converges to a stable equilibrium where both the healthy honeybee and the virus infected honeybees can persist (see the first row of Figure \ref{fig3:comparison1}). \\\\
 If there is no virus, the virus-free system \eqref{HM} (i.e., the parasitism dynamics of the healthy honeybee and the healthy mites) converges to the extinction of both species through  \emph{catastrophic event} (see the second row of Figure \ref{fig2:comparison1}).\\\\
 However, if honeybees, mites, and virus are all presented in the system (i.e., the full system \eqref{newHBD}), then both honeybee and mites go extinct (see Figure \ref{fig2:comparison2}).\\\\
 This implies that the presence of the unstable mite population can lead to the extinction of honeybees. \\
 \end{enumerate}
\begin{figure}[ht]
\centering
   \includegraphics[scale =.4] {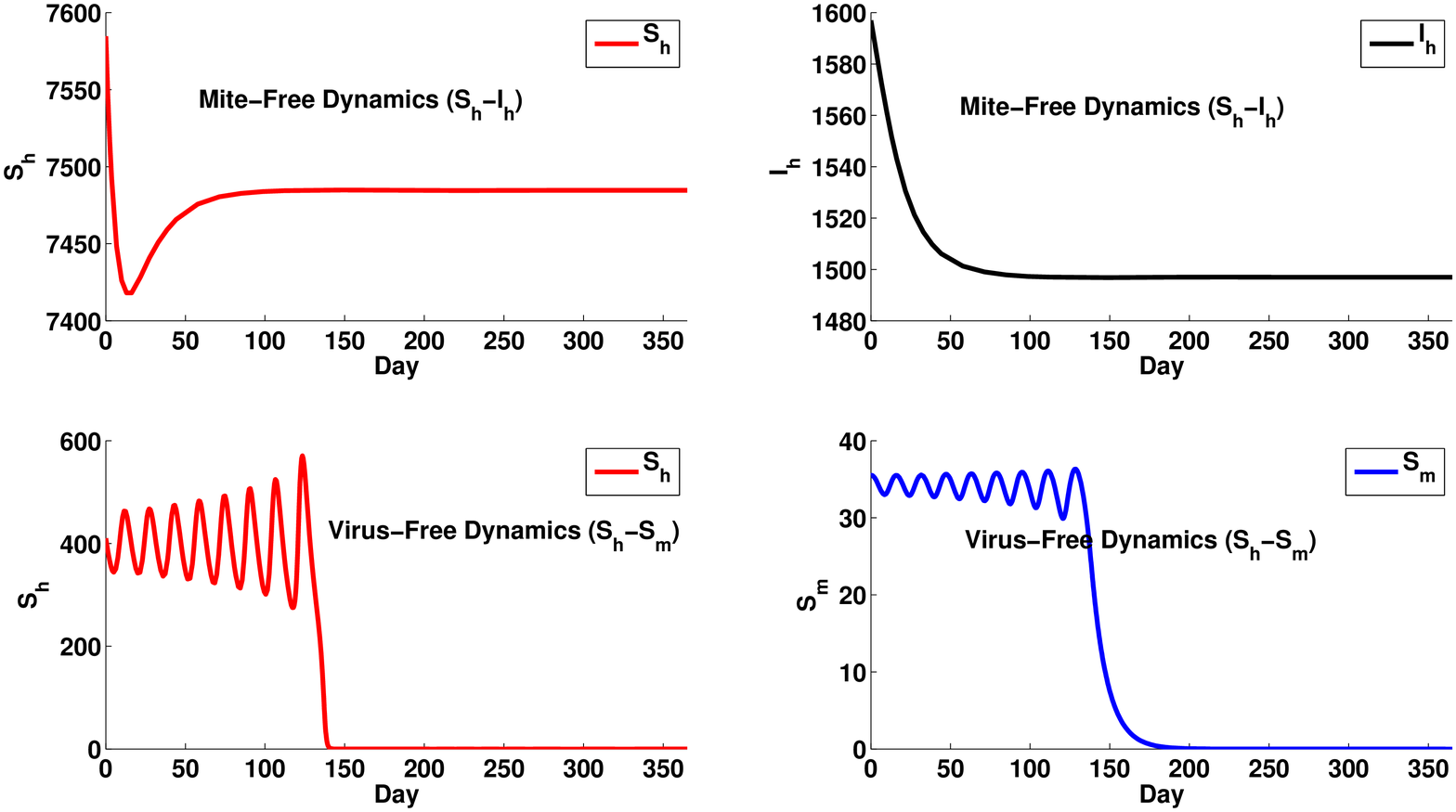}
  \caption{Population dynamics of the subsystems of the the honeybee-mite-virus model \eqref{newHBD} when $r = 1500; \rho=0.9; \hat{K}= 1600001; d_h =0.15; \mu_h=0.1; \alpha=0.05; c = 0.005; d_m=0.1;\mu_m=0.01; \beta_h=0.3; \hat{\beta}_{mh}=0.08; \tilde{\beta}_{mh}=0.001;\hat{\beta}_{hm}=0.03$. The left figure in the first row is the healthy honeybee population (the red curve) and the right figure in the first row is the virus infected honeybee population (the black curve) in the mite-free subsystem \eqref{HM} when  $S_h(0)=7684, I_h(0)=1700$.  The left figure in the second row is the healthy honeybee population (the red curve) and the right figure in the second row is the healthy mite population (the blue curve) in the virus-free subsystem \eqref{HV} when $S_h(0)=410,S_m(0)=35$.}\label{fig3:comparison1}
\end{figure}
\begin{figure}[ht]
\centering
 \includegraphics[scale =.4] {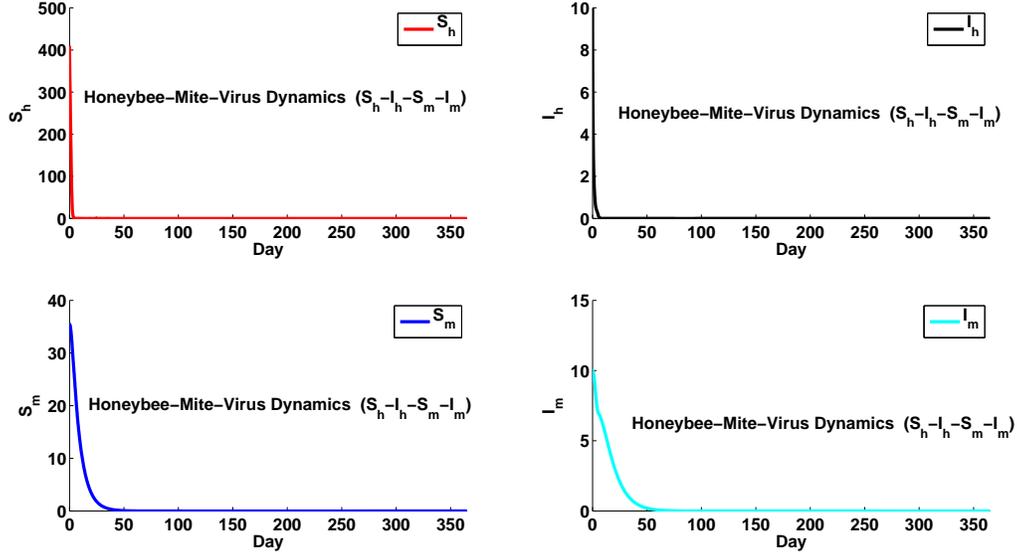}
\caption{Population dynamics of the honeybee-mite-virus model \eqref{newHBD} when $r=1500; K=160001; \rho=0.9; d_h=.15; \mu_h=0.1; \alpha=0.05; d_m=0.1; \mu_m=0.01;  c=0.005; \beta_h=.3; \hat{\beta}_{mh} = 0.08; \tilde{\beta}_{mh} = .001; \hat{\beta}_{hm} = 0.03$ and $S_h(0)=410, I_h(0)=10, S_m(0)=35, I_m(0)=10.$
  The healthy honeybee population $S_h$ is in red; the \modifyb{virus infected} honeybee population $I_h$ is in black; the healthy mite population $S_m$ is in blue; and the \modifyb{virus infected} mite population $I_m$ is in cyan.}\label{fig3:comparison2}
\end{figure}
\item Linking to CCD: In the absence of virus infections, the subsystem \eqref{HM} goes through the \emph{catastrophic event} which causes the extinction of honeybees (see the honeybee population in the red curve of Figure \ref{fig1:DF}). This property has been \modifyb{inherited by} the full system \eqref{newHBD} as honeybee population \modifyb{goes extinct} suddenly  for most initial conditions (see the honeybee population in the red curve of the fist left figure in \ref{fig3:comparison2}).\\
\end{enumerate}}
\modifyb{Our analysis and simulations suggest that it is important to include both mites and virus in studying the population dynamics of honeybees as we proposed the full system \eqref{newHBD} due to the synergistic effects of parasitism induced by mites and the virus infections as well as the \emph{catastrophic event} from the parasitism interactions between mites and honeybees. In other words, if we only consider the honeybee versus virus dynamics as described by the subsystem \eqref{HV}, or only consider the parasitism interactions between mites and honeybees as described by the subsystem \eqref{HM}, we are not able to capture the full mechanics that can lead to the extinction of honeybees or the biological implications on the persistence of disease in honeybees. }
In addition, the full system \eqref{newHBD} has rich dynamics which can \modifyb{possess} multiple attractors. Thus, depending on initial conditions, the full system \eqref{newHBD} can either experience extinction or have coexistence of  both honeybees and mites.  More sophisticated mathematical analysis is needed in order to understand the detailed dynamics. \\


\section{Discussion}

The association of virus infection with Varroa mites infestation in honeybee colonies causes great concern for researchers and beekeepers \cite{Chen:2006ab}. Many studies have suggested that Varroa mite infestations could be a key explanatory factor for the widespread increase in annual honeybee colony mortality and has been implicated as a contributing factor leading to CCD \cite{Mondet:2014aa}. In this paper, we derive and study a honeybee-virus-mite model by using the susceptible-infectious (SI) disease framework. Our proposed model includes the parasitic interaction between Varroa mites and honeybees as well as different virus transmission modes occurring at different phases of Varroa mites. More specifically, we use frequency-dependent transmission functions to model horizontal virus transmissions between Varroa mites and honeybees during the phoretic phase of mites and use a Holling Type I functional response to model parasitic interactions between Varroa mites and honeybees. We also use
a density-dependent transmission function to model the vertical virus transmission between Varroa mites and honeybees during the reproductive phase of Varroa mites. We also apply the frequency-dependent transmission function to model the horizontal virus transmission among honeybees.\\

\modifyb{\textbf{Summarize the main findings:} Our analytical and numerical results of the full system suggest the following:\\
\begin{enumerate}
\item Initial honeybee populations play an important role in its persistence  since its dynamics exhibits strong Allee effect in the absence of both parasites and virus (see Theorem \ref{th1:basics}). In addition, patterns of population dynamics are sensitive to initial conditions as suggested by our numerical results, i.e., depending on initial conditions, the full system can experience the \emph{catastrophic event} where honeybees collapse dramatically, or both mites and honeybees can coexist but exhibit different dynamical patterns.\\
\item In the absence of Varroa mites, the honeybee and disease system has only equilibrium dynamics (see Theorem \ref{th3:MF} and the first row of Figure \ref{fig2:comparison1}-\ref{fig2:comparison2}). In the presence of Varroa mites, the synergistic effects make disease and mites more persistent (see Theorem \ref{th4:HVMI}-Theorem \ref{th6:global}), and thus, difficult to control (see Theorem \ref{th2:per-ext} and the related arguments). In addition, the synergistic effects make disease and mites  can also drive both honeybees and mites go extinct (see Figure \ref{fig3:comparison2}).\\ 
\item In the absence of virus infections, the Varroa mite and honeybee system can be destabilized by the low adult workers to brood ratio in the colony where the system has oscillating dynamics leading to a sudden extinction of all species (see Theorem \ref{th5:DF} and Figure  \ref{fig1:DF}). This dynamical property is called the \emph{catastrophic event} which has been \modifyb{inherited by} the full system when the virus infection is present (see Figure \ref{fig3:comparison2}). This phenomenon could be linked to CCD which has been observed in honeybee colonies where low hive to brood ratio could be a contributing factor.\\
\item Our numerical simulations suggest that large values of virus transmission rates can drive the extinction of healthy mites while extremely large values can lead to the extinction of all species.\\
\end{enumerate}}

\textbf{Identify the contribution to the broader field}: Our theoretical results in Theorem \ref{th1:basics} imply that the persistence of honeybees requires \modifyr{a critical threshold of the adult workers to brood ratio. This illustrates previous observations that a proper adult worker to brood ratio is needed to have a successful honeybee colony \cite{Schmickl:2004aa}.} Our current work indeed \modifyb{provides} many useful insights on our proposed questions regarding the complicated synergistic effects of  Varroa mites and the associated virus infections on the honeybee population dynamics and the persistence as well as how these synergistic effects may potentially contribute to CCD of honeybee colonies. \modifyg{Need more help from Gloria}\\

 \textbf{Compare the findings with other works:} \modifyc{As we mentioned in our introduction, the work of Sumpter and Martin \cite{Sumpter:2004aa} and Ratti \textit{et al.} \cite{Ratti:2012aa} are most relevant modeling papers for our study proposes. However, Sumpter and Martin \cite{Sumpter:2004aa} assumed that the mites' population is constant and the virus transmission occurs only through Varroa mites.  This assumption prevents us to study how mites' population affect the virus transmission and the honeybee population dynamics. The model of Ratti \textit{et al.} \cite{Ratti:2012aa} does not include the virus transmission through honeybees themselves which leads to the extinction of virus when Varroa mites go extinction. This is not realistic.  Our honeybee-mite-virus model allows us to investigate the honeybee population dynamics in the absence of both the Varroa mites and virus (the honeybee dynamics \eqref{H}),  in the absence of the virus (the virus-free dynamics \eqref{HM}), in the absence of the Vorroa mites (the mite-free dynamics \eqref{HV}, where the virus is transmitted through other contacts), and in the presence of both Varroa mites and virus to explore the synergistic effects of  Varroa mites and the associated virus infections on the honeybee population dynamics as well as the potential mechanisms contributing to CCD of honeybee colonies. More specifically, our model is able to provide sufficient conditions that can guarantee the persistence of honeybee population and conditions that lead to the extinction of honeybees.} \\
 
 \modifyr{We would like to mention the recent work by Perry \textit{et al.} \cite{Perry:2015aa} which examined the social dynamics underlying the dramatic colony failure with an aid of a honeybee population model. Their model does not include population dynamics of Vorroa mites or virus infections, but it does includes bee foraging performance varying with age, and displays dynamics of colony population collapse that are similar to field reports of CCD. We plan to expand our current modeling framework to include the detailed stage-structure model of honeybees with social dynamics. } \\
\modifyg{Even if nutrients are not limiting, as assumed in the VARROAPOP model, colonies go to extinction due to reduced longevity of adult workers, adult population decline, and a reduced ratio of adults to brood. This ultimately reduces brood rearing so that bees die at rates that are greater than emergence of new adults. This begins the downward spiral in the adult population that ultimately leads to colony death.}\\


\textbf{Suggest options for further research:} \modifyb{In our proposed model, we assume that a fixed ratio of adult workers to brood which does not reflect the reality since this ratio varies with the availability of nutrient that depends on the quality and availability of pollen. The mortality of honey bees is resulted from different risk factors such as parasites, pathogen, viruses, pesticides, nutrition and environmental changes \cite{Speybroeck:2010aa}.  Among these factors, honey bee nutrition is one of the most fundamental factors which impacts honey bee health and influences the capabilities of honey bees to combat different stressors. Nutrients such as protein from pollen is essential in fighting parasites and diseases, and in maintaining the high adult worker to brood ratio in honeybees. This is because that honeybees solely depend on pollen and honey for satisfying their daily needs where pollen provides protein, essential amino acids, vitamins and mineral for bees while honey or nectar are their carbohydrate resources \cite{Atkins:1975aa}. Cage studies showed that honey bees could survive for a long time without pollen provision \cite{Barker:1974aa}, but pollen feeding significantly prolonged their lifespan \cite{Manning:2007aa,Maurizio:1954aa,Schmidt:1995aa,Schmidt:1987aa}.}  \\

\modifyb{In addition, intensive research has investigated the connection between nutrition and honey bee disease and stress resistance. Both cage studies and field studies indicated that bees with poor nutrition were under more stress (Wang et al. unpublished data) \cite{Wang:2014aa}, more susceptible to Nosema and Varroa destructor, and had shorter lifespan \cite{Eischen:2008aa,Rinderer:1977aa}. DeGrandi-Hoffman found a similar result that bees with good pollen or protein nutrition were more resistant to many different viruses \cite{DeGrandi-Hoffman:2010aa} . Furthermore, studies on molecular mechanisms suggested that pollen nutrition may positively affect antimicrobial peptides and improve immune defensive response to parasites \cite{Di-Pasquale:2013aa}. Therefore, a more realistic model that includes nutrient/brood dynamics, explicit division of labor, and seasonal effects, are needed. This is our on-going project.} \\

\section{Proof}

\subsection*{Proof of Theorem \ref{th1:basics}}
\begin{proof}According to Theorem A.4 (p.423) of Thieme (2003), we can conclude that Model \eqref{newHBD} is positive invariant in $X$.
Let $d=\min\{d_h, d_m\} ,\,N_h=S_h+ I_h,\, N_m=S_m+ I_m,\, \mbox{ and } N=cN_h + N_m$. Then we have
$$\begin{array}{lcl}
N'=cN_h' + N_m'&=&\frac{rc\left(S_h+\rho I_h\right)^2}{\hat{K}+\left(S_h+\rho I_h\right)^2}-cd_h N_h-c\mu_h I_h- d_m N_m-\mu_m I_m\\\\
&\leq&\frac{rc(N/c)^2}{\hat{K}+(N/c)^2}-\min\{d_h, d_m\} (cN_h + N_m)=\frac{rcN^2}{c^2\hat{K}+N^2}-\min\{d_h, d_m\} N=\frac{rcN^2}{c^2\hat{K}+N^2}-d N
\end{array}$$ which implies that $\limsup_{t\rightarrow\infty} N(t) \leq c\frac{\frac{r}{d}+\sqrt{\left(\frac{r}{d}\right)^2-4\hat{K}}}{2}=N^*$ with implication that $\limsup_{t\rightarrow\infty} N(t)=0$ if either $\frac{r}{d}<2\sqrt{\hat{K}}$ or $N(0)<N^c=c\frac{\frac{r}{d}+\sqrt{\left(\frac{r}{d}\right)^2-4\hat{K}}}{2}$ holds. 
The arguments above also imply that if $\frac{r}{d}>2\sqrt{\hat{K}}$, then we have
$$\limsup_{t\rightarrow\infty} N_h(t)\leq N^*/c \mbox{ and }\limsup_{t\rightarrow\infty} N_m(t)\leq N^*.$$
Similarly, we have follows for $N_h$:
$$\begin{array}{lcl}
N_h'=S_h'+ I_h'&=&\frac{r\left(S_h+\rho I_h\right)^2}{\hat{K}+\left(S_h+\rho I_h\right)^2}-d_h N_h-\mu_h I_h- \alpha N_h N_m\\\\
&\leq&\frac{rN_h^2}{\hat{K}+N_h^2}-d_h N_h
\end{array}$$ which implies that $$\limsup_{t\rightarrow\infty} N_h(t)\leq \frac{\frac{r}{d_h}+\sqrt{\left(\frac{r}{d_h}\right)^2-4\hat{K}}}{2}=\bar{N}_h^* \mbox{ if } \frac{r}{d_h}>2\sqrt{\hat{K}}$$ and 
$$\limsup_{t\rightarrow\infty} N_h(t)=0 \mbox{ if } \frac{r}{d_h}<2\sqrt{\hat{K}} \mbox{ or } N_h(0)<S^c_h=\frac{\frac{r}{d_h}-\sqrt{\left(\frac{r}{d_h}\right)^2-4\hat{K}}}{2}.$$
On the other hand, we have
$$\begin{array}{lcl}
N_h'=S_h'+ I_h'&\geq &\frac{r\rho^2\left(S_h/\rho+ I_h\right)^2}{\hat{K}+\rho^2\left(S_h/\rho+ I_h\right)^2}-(d_h +\mu_h+\alpha N_m)N_h\\\\
&\geq&\frac{r\rho^2N_h^2}{\hat{K}+\rho^2N_h^2}-(d_h +\mu_h+\alpha N^*)N_h\\\\
&\geq&\frac{rN_h^2}{\hat{K}/\rho^2+N_h^2}-(d_h +\mu_h+\alpha N^*)N_h
\end{array}.$$
Therefore, apply the comparison theorem, we can conclude that
$$\liminf_{t\rightarrow\infty} N_h(t)\geq \frac{\frac{r}{d_h +\mu_h+\alpha N^*}+\sqrt{\left(\frac{r}{d_h +\mu_h+\alpha N^*}\right)^2-4\hat{K}/\rho^2}}{2}=\underline{N}^*_h $$\mbox{ if } the following inequalities hold
$$ \frac{r}{d_h +\mu_h+\alpha N^*}>2\sqrt{\frac{\hat{K}}{\rho^2}}\mbox{ and } N_h(0)>\underline{N}^c_h=\frac{\frac{r}{d_h +\mu_h+\alpha N^*}-\sqrt{\left(\frac{r}{d_h +\mu_h+\alpha N^*}\right)^2-4\hat{K}/\rho^2}}{2}.$$
The discussions above  provide sufficient conditions that allow $N_h$ being persistent, i.e., when the inequalities $\frac{r}{d_h +\mu_h+\alpha N^*}>2\sqrt{\frac{\hat{K}}{\rho^2}}\mbox{ and } N_h(0)>\underline{N}^c_h$ hold. This implies that $S_h$ is also persistent under this condition since all species go extinct if $S_h$ goes extinct. However, this persistence condition does not provide an estimate of $S_h$. To explore an estimation of $S_h$, we look at the population of $S_h$.\\


Recall that if $\frac{r}{d_h +\mu_h+\alpha N^*}>2\sqrt{\frac{\hat{K}}{\rho^2}}$ and $N_h(0)>\underline{N}^c_h$, then we have the following inequalities
$$\underline{N}^*_h\leq \liminf_{t\rightarrow\infty} N_h(t)\leq \liminf_{t\rightarrow\infty} N(t)/c\leq \limsup_{t\rightarrow\infty} N(t)/c\leq N^*/c.$$
If the inequalities $\frac{r}{d_h+\beta_h+\modifyr{\frac{\hat{\beta}_{mh} N^*}{\underline{N}^*_h}}+(\tilde{\beta}_{mh}+\alpha)(N^*-c\underline{N}^*_h)}>2\sqrt{\hat{K}}$ and $$S_h(0)>\frac{\frac{r}{d_h+\beta_h+\modifyr{\frac{\hat{\beta}_{mh} N^*}{\underline{N}^*_h}}+(\tilde{\beta}_{mh}+\alpha)(N^*-c\underline{N}^*_h)}-\sqrt{\left(\frac{r}{d_h+\beta_h+\modifyr{\frac{\hat{\beta}_{mh} N^*}{\underline{N}^*_h}}+(\tilde{\beta}_{mh}+\alpha)(N^*-c\underline{N}^*_h)}\right)^2-4\hat{K}}}{2}$$ hold, then we have following inequalities
$$\begin{array}{lcl}
S_h'&=& \frac{r\left(S_h+\rho I_h\right)^2}{\hat{K}+\left(S_h+\rho I_h\right)^2}-d_h S_h-\frac{\beta_{h} S_hI_h}{S_h+I_h}-\modifyr{\frac{\hat{\beta}_{mh} S_hI_m}{S_h+I_h}}-\tilde{\beta}_{mh}S_hI_m-\alpha S_h (S_m+I_m)\\\\
&>&\frac{r\left(S_h\right)^2}{\hat{K}+\left(S_h\right)^2}-d_h S_h-\beta_{h} S_h-\modifyr{\frac{\hat{\beta}_{mh} I_m}{N_h}}S_h-\tilde{\beta}_{mh}S_hN_m-\alpha S_h N_m\\\\
&>&\frac{r\left(S_h\right)^2}{\hat{K}+\left(S_h\right)^2}-\left[d_h+\beta+\modifyr{\frac{\hat{\beta}_{mh} N^*}{\underline{N}^*_h}}+(\tilde{\beta}_{mh}+\alpha)N_m\right]S_h\\\\
&=&\frac{r\left(S_h\right)^2}{\hat{K}+\left(S_h\right)^2}-\left[d_h+\beta+\modifyr{\frac{\hat{\beta}_{mh} N^*}{\underline{N}^*_h}}+(\tilde{\beta}_{mh}+\alpha)(N-cN_h)\right]S_h\\\\
&>&\frac{r\left(S_h\right)^2}{\hat{K}+\left(S_h\right)^2}-\left[d_h+\beta+\modifyr{\frac{\hat{\beta}_{mh} N^*}{\underline{N}^*_h}}+(\tilde{\beta}_{mh}+\alpha)(N^*-c\underline{N}^*_h)\right]S_h
\end{array}$$which implies that
$$\liminf_{t\rightarrow\infty}S_h(t)\geq S_h^*=\frac{\frac{r}{d_h+\beta_h+\modifyr{\frac{\hat{\beta}_{mh} N^*}{\underline{N}^*_h}}+(\tilde{\beta}_{mh}+\alpha)(N^*-c\underline{N}^*_h)}+\sqrt{\left(\frac{r}{d_h+\beta_h+\modifyr{\frac{\hat{\beta}_{mh} N^*}{\underline{N}^*_h}}+(\tilde{\beta}_{mh}+\alpha)(N^*-c\underline{N}^*_h)}\right)^2-4\hat{K}}}{2}.$$
Therefore, we can conclude that $S_h$ is persistent with the following properties:
$$S^*_h\leq \liminf_{t\rightarrow\infty} S_h(t)\leq \liminf_{t\rightarrow\infty} N_h(t)\leq \liminf_{t\rightarrow\infty} \frac{N(t)}{c}\leq \limsup_{t\rightarrow\infty} N(t)/c\leq N^*/c$$ if the following inequalities hold
\begin{enumerate}
\item $\frac{r}{2\sqrt{\hat{K}}}>\max\Big\{d_h+\beta_h+\modifyr{\frac{\hat{\beta}_{mh} N^*}{\underline{N}^*_h}}+(\tilde{\beta}_{mh}+\alpha)(N^*-\underline{N}^*_h),\, \frac{d_h +\mu_h+\alpha N^*}{\rho}\Big\}$ \mbox{ with } $ N_h(0)\geq S_h(0)>S_h^c $.
\end{enumerate}
It is easy to check that the extinction equilibrium $E_0=(0,0,0,0)$ is always local stable. We omit the details. Based on the discussion of the upper bound of the total population $N$ and the total population of honeybees $N_h$, we can conclude that the system \eqref{newHBD} converges to $E_0$ globally if $2\sqrt{\hat{K}}>\frac{r}{d_h}$ 
holds; while if the initial population satisfies either $N(0)<N^c$ or $N_h(0)<\bar{N}^c_h$, then the system \eqref{newHBD} converges to $E_0$ locally.

\end{proof}

\subsection*{Proof of Theorem \ref{th2:per-ext}}
\begin{proof}
Now we consider the population of $N_m$. From Model \eqref{newHBD}, we obtain the following inequalities:
$$\begin{array}{lcl}
c\alpha N_h N_m-(d_m+\mu_m) N_m \leq N_m'=S_m'+ I_m'&=&c\alpha N_h N_m-d_m N_m-\mu_m I_m \leq N_m\left(c\alpha N_h -d_m \right)
\end{array}$$ which implies that if $N^*<\frac{d_m}{\alpha }$, then we have 
$$\begin{array}{lcl}
 N_m'&\leq& N_m\left[c\alpha N_h -d_m \right]=N_m\left[c\alpha (N-N_m)/c -d_m \right]\\
 &\leq& N_m\left[\alpha N^*-d_m-\alpha N_m \right]<N_m\left[\alpha N^*-d_m \right]
\end{array}.$$ This indicates that $\limsup_{t\rightarrow\infty} N_m(t)=0.$\\


 Assume that the following inequalities hold
$\frac{r}{d_h}>2\sqrt{\hat{K}},\,  \mbox{ and }S_h(0)>{S}_h^c.$ Then according to Theorem \ref{th1:basics}, we have 
$$\limsup_{t\rightarrow\infty}S_h(t)\leq \limsup_{t\rightarrow\infty}N_h(t)\leq \bar{N}_h^*.$$
Now let us focus on the population of $S_m$. Notice that we have the following inequalities when time is large enough,
$$\begin{array}{lcl}
{S_m'}&=&{S_m}\left[c\alpha S_h-\frac{\hat{\beta}_{hm}I_h}{S_h+I_h}-d_m\right]\leq S_m\left[c\alpha N_h-d_m\right]\leq S_m\left[c\alpha \bar{N}^*_h-d_m\right]
 \end{array}.$$
 This implies that if $\bar{N}_h^*<\frac{d_m}{c\alpha}$ and $\frac{r}{d_h}>2\sqrt{\hat{K}},\,  \mbox{ and }S_h(0)>{S}_h^c$, the healthy mite population $S_m$ goes extinct while the total population of honeybees persists, i.e., 
$$\limsup_{t\rightarrow\infty} S_m(t)= 0.$$

Now we look at the population dynamics of $I_m$ and $I_h$. Let $I=cI_h+I_m$. From Model \eqref{newHBD}, then we have the following equations:
$$\begin{array}{lcl}
I'&=&(cI_h+I_m)'=S_h\left[\frac{c\beta_{h}I_h}{S_h+I_h}+\modifyr{\frac{c\hat{\beta}_{mh} I_m}{S_h+I_h}}+c\tilde{\beta}_{mh}I_m\right]-c\alpha I_h (S_m+I_m)-c(d_h+\mu_h)I_h\\
&&+c\alpha\left[I_h (S_m+I_m)+S_hI_m\right]+\frac{\hat{\beta}_{hm}I_hS_m}{S_h+I_h}-(d_m +\mu_m)I_m\\\\
&=&I_h\left[\frac{c\beta_{h}S_h}{N_h}+\frac{c\hat{\beta}_{hm}S_m}{N_h}\right]+I_m\left[\modifyr{\frac{c\hat{\beta}_{mh} S_h}{N_h}}+c\tilde{\beta}_{mh}+c\alpha S_h\right]-c(d_h+\mu_h)I_h-(d_m +\mu_m)I_m
\end{array}.$$
This implies that
$$I'\geq I\left[ \min\Big\{\frac{\beta_{h}S_h}{N_h}+\frac{\hat{\beta}_{hm}S_m}{N_h},\modifyr{\frac{c\hat{\beta}_{mh} S_h}{N_h}}+c\tilde{\beta}_{mh}+c\alpha S_h\Big\}-\max\Big\{(d_h+\mu_h),(d_m +\mu_m)\Big\} \right].$$ 
Assume that $\frac{r}{2\sqrt{\hat{K}}}>d_h+\beta_h+\modifyr{\frac{\hat{\beta}_{mh} N^*}{\underline{N}^*_h}}+(\tilde{\beta}_{mh}+\alpha)(N^*-\underline{N}^*_h)\mbox{ and }S_h(0)>{S}_h^c,$ then according to Theorem \ref{th1:basics}, we have 
$$\liminf_{t\rightarrow\infty} S_h(t)\geq {S}^*_h.$$ 
Therefore, if $\frac{\min\Big\{\beta_{h},\modifyr{c\hat{\beta}_{mh} }+c\tilde{\beta}_{mh}+c\alpha {S}^*_h\Big\}}{\max\Big\{(d_h+\mu_h),(d_m +\mu_m)\Big\}}\geq1$, then we have
{\small$$\begin{array}{lcl}
\frac{I'}{I}\Big\vert_{I_h=I_m=0}&\geq&\left[ \min\Big\{\frac{\beta_{h}S_h}{N_h},\modifyr{c\hat{\beta}_{mh}}+c\tilde{\beta}_{mh}+c\alpha S_h\Big\}-\max\Big\{d_h+\mu_h,(d_m +\mu_m)\Big\} \right]\Big\vert_{I_h=I_m=0}\\\\
&>&\left[ \min\Big\{\beta_{h},\modifyr{c\hat{\beta}_{mh}}+c\tilde{\beta}_{mh}+c\alpha {S}^*_h\Big\}-\max\Big\{(d_h+\mu_h),(d_m +\mu_m)\Big\} \right]\geq0
\end{array}$$}which implies that the disease $I$ persists by applying the average Lynapunov theorem (Hutson 1984). \\

On the other hand, we have the following inequalities:
{\small$$\begin{array}{lcl}
I'&\leq&I\left[ \max\Big\{\frac{\beta_{h}S_h}{N_h}+\frac{\hat{\beta}_{hm}S_m}{N_h},\frac{c\hat{\beta}_{mh}S_h}{N_h}+c\tilde{\beta}_{mh}+c\alpha S_h\Big\}-\min\Big\{(d_h+\mu_h),(d_m +\mu_m)\Big\} \right]\\\\
\frac{I'}{I}\Big\vert_{I_h=I_m=0}&<&\left[ \max\Big\{\beta_{h}+\frac{\hat{\beta}_{hm}S_m}{N_h},\frac{c\hat{\beta}_{mh}S_h}{N_h}+c\tilde{\beta}_{mh}+c\alpha S_h\Big\}-\min\Big\{(d_h+\mu_h),(d_m +\mu_m)\Big\} \right]
\end{array}.$$}
Assume that $\frac{r}{2\sqrt{\hat{K}}}>\frac{d_h+\mu_h+\alpha N^*}{\rho},\,\,N_h(0)>\underline{N}^c_h$. Then according to Theorem \ref{th1:basics}, we have the following inequalities:
$$\underline{N}^*_h\leq\liminf_{t\rightarrow\infty} N_h(t)\leq \limsup_{t\rightarrow\infty} N_h(t)\leq \bar{N}^*_h\leq\limsup_{t\rightarrow\infty} N(t)/c=N^*/c.$$
This implies the following inequality when $\frac{\max\Big\{\beta_{h}+\frac{\hat{\beta}_{hm}N^*}{\underline{N}^*_h},\modifyr{c\hat{\beta}_{mh}}+c\tilde{\beta}_{mh}+c\alpha \bar{N}^*_h\Big\}}{\min\Big\{(d_h+\mu_h),(d_m +\mu_m)\Big\}}<1$, then we have
{\small$$\begin{array}{lcl}
\frac{I'}{I}\Big\vert_{I_h=I_m=0}&<&\left[ \max\Big\{\beta_{h}+\frac{\hat{\beta}_{hm}N^*}{\underline{N}^*_h},\modifyr{c\hat{\beta}_{mh}}+c\tilde{\beta}_{mh}+c\alpha \bar{N}^*_h\Big\}-\min\Big\{(d_h+\mu_h),(d_m +\mu_m)\Big\} \right]<0
\end{array}$$} which implies that the disease goes extinct, i.e.,
$$\lim_{t\rightarrow\infty} I(t)=0.$$

\end{proof}
\subsection*{Proof of Theorem \ref{th5:DF}}
\begin{proof}If $M=0$, the virus-free subsystem \eqref{HM} reduces to the only healthy honeybee population:
$$S_h'= \frac{rS_h^2}{\hat{K}+S_h^2}-d_h S_h$$ which leads to the following three boundary equilibria if $\frac{r}{2\hat{K}}>d_h$:
$$(0,0),\,\,(\bar{N}^c_h,0),\,\, \mbox{ and } (\bar{N}^*_h,0).$$
The local stability can be easily determined by the eigenvalues evaluated at its Jacobian matrix. Simple algebraic calculations show that $(\bar{N}^c_h,0)$ is a saddle if $\bar{N}^c_h<\frac{d_m}{\alpha c}$ while it is a source if $\bar{N}^c_h>\frac{d_m}{\alpha c}$; and $(\bar{N}^*_h,0)$ is a sink if $\bar{N}^*_h<\frac{d_m}{\alpha c}$ while it is a saddle if $\bar{N}^*_h>\frac{d_m}{\alpha c}$.\\

Now let $(H,M)$ be an interior equilibrium of the virus-free subsystem \eqref{HM}, then we have the following equations hold:
$$\begin{array}{lcl}
0&=&c\alpha HM-d_m M \Leftrightarrow H=\frac{d_m}{c\alpha}\\\\
0&=& \frac{rH^2}{\hat{K}+H^2}-d_h H-\alpha HM \Leftrightarrow M=\frac{1}{\alpha}\left(\frac{rH}{\hat{K}+H^2}-d_h\right)
\end{array}.$$
which gives the unique interior equilibrium $(H^*, M^*)=\left(\frac{d_m}{c\alpha},\frac{1}{\alpha}\left(\frac{rH^*}{\hat{K}+(H^*)^2}-d_h\right)\right)$ provided that $\bar{N}^c_h<\frac{d_m}{\alpha c}<\bar{N}^*_h$.\\

The local stability of $(H^*, M^*)$ is determined by the eigenvalues $\lambda_i,i=1,2$ of the following Jacobian matrix of \eqref{HM}:
\bae\label{JHM}
J_{HM}&=&\left(\begin{array}{cc}
\frac{rH^*\left(\hat{K}-(H^*)^2\right)}{\left(\hat{K}+(H^*)^2\right)^2}&-\alpha H^*\\\\
\alpha c M^*&0\end{array}\right)
\eae which gives the following two equations:
\bae\label{tJHM}
\begin{array}{lcl}
\lambda_1+\lambda_2&=&\frac{rH^*\left(\hat{K}-(H^*)^2\right)}{\left(\hat{K}+(H^*)^2\right)^2}\\\\
\lambda_1\lambda_2&=&c\alpha^2 H^*M^*
\end{array}.
\eae
Therefore, we can conclude that $(H^*, M^*)$ is a sink if  $H^*>\sqrt{\hat{K}}$ while $(H^*, M^*)$ is a source if  $H^*<\sqrt{\hat{K}}$.\\ 
The discussion above implies follows:
\begin{enumerate}
\item If $\frac{r}{2\sqrt{\hat{K}}}<d_h$, then the extinction equilibrium $(0,0)$ is the only locally stable equilibrium for the virus-free subsystem \eqref{HM}. Thus it is globally stable.
\item If $\frac{r}{2\sqrt{\hat{K}}}>d_h$ and either $\bar{N}^c_h>\frac{d_m}{\alpha c}$ or $\frac{d_m}{\alpha c}>\bar{N}^*_h$, then the virus-free subsystem \eqref{HM} has three boundary equilibria: $(0,0),\,(S^c_h,0),\,(\bar{N}_h^*,0)$ where $(0,0)$ is always locally stable.
\begin{itemize}
\item If $\bar{N}^c_h>\frac{d_m}{\alpha c}$, then $(S^c_h,0)$ is a source and $(\bar{N}_h^*,0)$ is a saddle. Since \eqref{HM} is a two-D ode system, it has global stability at the extinction equilibrium $(0,0)$ according to the Poincare-Bendison theorem \cite{Perko:2001aa}.
\item If $\bar{N}^*_h<\frac{d_m}{\alpha c}$, then $(S^c_h,0)$ is a saddle and $(\bar{N}_h^*,0)$ is a sink. This implies that \eqref{HM} has two locally asymptotically stable boundary equilibria $(0,0)$ and $(\bar{N}_h^*,0)$ which are reserved as the only two attractors for the model.\end{itemize}
\item If $\frac{r}{2\sqrt{\hat{K}}}>d_h$ and $\bar{N}^c_h<\frac{d_m}{\alpha c}<\bar{N}^*_h$, then  \eqref{HM} has three boundary equilibria $(0,0),\,(S^c_h,0),\,(\bar{N}_h^*,0)$ and the unique interior equilibrium $(H^*,M^*)$ where $(0,0)$ is locally stable, both $(S^c_h,0)$ and $(\bar{N}_h^*,0)$ are saddle nodes. The local stability of $(H^*,M^*)$ is determined by the sign of $H^*-\sqrt{\hat{K}}$: if $H^*-\sqrt{\hat{K}}>0$, then $(H^*,M^*)$ is locally asymptotically stable while if $H^*-\sqrt{\hat{K}}<0$, then $(H^*,M^*)$ is a source.
\end{enumerate}
Define $\mathcal R_0^M=\frac{\bar{N}_h^*}{H^*}$. Then we can conclude that \eqref{HM}  has no interior equilibrium if $\mathcal R_0^M<1$ or $\mathcal R_0^M>\frac{\bar{N}_h^*}{S^c_h}>1$; and  \eqref{HM}  has a unique interior equilibrium $(H^*,M^*)$ if $1<\mathcal R_0^M<\frac{\bar{N}_h^*}{S^c_h}.$ More specifically, $(H^*,M^*)$ is locally stable if $1<\mathcal R_0^M<\frac{\bar{N}_h^*}{\hat{K}}$ while it is a source if $\mathcal R_0^M>\max\Big\{1, \frac{\bar{N}_h^*}{\hat{K}}\Big\}$.
\end{proof}
\subsection*{Proof of Theorem \ref{th3:MF}}
\begin{proof}If $I=0$, the mite-free subsystem \eqref{HV} reduces to the only healthy honeybee population:
$$S_h'= \frac{rS_h^2}{\hat{K}+S_h^2}-d_h S_h$$ which leads to the following three boundary equilibria if $\frac{r}{2\hat{K}}>d_h$:
$$(0,0),\,\,(\bar{N}^c_h,0),\,\, \mbox{ and } (\bar{N}^*_h,0).$$
The local stability can be easily determined by the eigenvalues evaluated at its Jacobian matrix. Define $\mathcal R_0^V=\frac{\beta_h}{d_h+\mu_h}$. The simple algebraic calculations show that if $\mathcal R_0^V<1$, $(\bar{N}^c_h,0)$ is a saddle and $(\bar{N}^*_h,0)$ is a sink; while if $\mathcal R_0^V>1$, then $(\bar{N}^c_h,0)$ is a source and $(\bar{N}^*_h,0)$ is a saddle. According to Theorem \ref{th1:basics}, the extinction equilibrium $(0,0)$ is always a sink.\\

Let $(S, I)$ be an interior equilibrium of the mite-free subsystem \eqref{HV}, then it satisfies the following equations:
\bae\label{HVeq}
\begin{array}{lcl}
\frac{\beta_{h}S}{S+I}-(d_h+\mu_h)=0&\Rightarrow &S=\frac{d_h+\mu_h}{\beta_{h}-(d_h+\mu_h)}I=\frac{I}{\frac{\beta_h}{d_h+\mu_h}-1}=aI\\\\
\frac{r (S+\rho I)^2}{\hat{K}+(S+\rho I)^2}-d_h S-S\frac{\beta_{h}I}{S+I}=0&\Rightarrow& \frac{r I}{\frac{\hat{K}}{(a+\rho)^2}+I^2}-a\left(d_h + \frac{\beta_{h}}{1+a}\right)=0\\\\

\frac{r (S+\rho I)^2}{\hat{K}+(S+\rho I)^2}-d_h S-S\frac{\beta_{h}I}{S+I}=0&\Rightarrow& \frac{r (S+\rho I)^2}{\hat{K}+(S+\rho I)^2}-d_h S-(d_h+\mu_h)I=0\\\\
&\Rightarrow &\frac{r I}{\frac{\hat{K}}{(a+\rho)^2}+I^2}-\left( (a+1)d_h +\mu_h\right)=\frac{r I}{\frac{\hat{K}}{(a+\rho)^2}+I^2}-\tilde{d}=0\\\\
&\Rightarrow &\frac{a\beta_{h}}{1+a}=d_h+\mu_h.
\end{array}
\eae
Therefore, if $\mathcal R_0^V=\frac{\beta_h}{d_h+\mu_h}>1$ and $\frac{r}{\tilde{d}}>\frac{2\sqrt{\hat{K}}}{a+\rho}$, then we can conclude that the mite-free subsystem \eqref{HV} can have two interior equilibria $(S_h^k,I_h^k), k=1,2$ where $$a=\frac{1}{\frac{\beta_h}{d_h+\mu_h}-1},\,\,d=(a+1) d_h+\mu_h=d_h\left(\frac{\frac{\beta_h}{d_h+\mu_h}}{\frac{\beta_h}{d_h+\mu_h}-1}\right)+\mu_h,$$\, and $$I_h^{1}=\frac{\frac{r}{\tilde{d}}-\sqrt{\left(\frac{r}{\tilde{d}}\right)^2-4\frac{\hat{K}}{(a+\rho)^2}}}{2},\,\,I_h^{2}=\frac{\frac{r}{\tilde{d}}+\sqrt{\left(\frac{r}{\tilde{d}}\right)^2-4\frac{\hat{K}}{(a+\rho)^2}}}{2},\,\, S_h^k= a I_h^k, k=1,2.$$
Now we examine the local stability of these two interior equilibria provided that they exist. The Jacobian matrix of \eqref{HV} evaluated at the interior equilibrium $(S,I)=(aI,I)$ can be expressed as follows:
\bae\label{JHV}
J_{HV}&=&\left(\begin{array}{cc}
\frac{2r\hat{K}(a+\rho)I}{((a+\rho)^2I^2+\hat{K})^2}-d_h-\frac{\beta_h}{(1+a)^2}&\frac{2r\rho \hat{K}(a+\rho)I}{((a+\rho)^2I^2+\hat{K})^2}-\frac{a^2\beta_h}{(1+a)^2}\\\\
\frac{\beta_h}{(1+a)^2}&-\frac{a\beta_h}{(1+a)^2}
\end{array}\right)
\eae whose two eigenvalues $\lambda_i,i=1,2$ satisfy the following equations:
\bae\label{tJHV}
\begin{array}{lcl}
\lambda_1+\lambda_2&=&\frac{2r\hat{K}(a+\rho)I}{((a+\rho)^2I^2+\hat{K})^2}-d_h-\frac{\beta_h}{1+a}
=\frac{2\hat{K}}{((a+\rho)^2I^2+\hat{K})(a+\rho)}\frac{r (a+\rho)^2I}{(a+\rho)^2I^2+\hat{K}}-d_h-\frac{\beta_h}{1+a}\\\\
&=&\frac{2\hat{K}}{((a+\rho)^2I^2+\hat{K})(a+\rho)} a \left(d_h+\frac{\beta_h}{1+a}\right)-d_h-\frac{\beta_h}{1+a}=\left(d_h+\frac{\beta_h}{1+a}\right)\left[\frac{2a\hat{K}}{((a+\rho)^2I^2+\hat{K})(a+\rho)} -1\right]\\\\
&<&\left(d_h+\frac{\beta_h}{1+a}\right)\left[\frac{2\hat{K}}{(a+\rho)^2I^2+\hat{K}} -1\right]
\end{array}
\eae and 
\bae\label{dJHV}
\begin{array}{lcl}
\lambda_1\lambda_2&=&\frac{\beta_h}{(1+a)^2}\left[-\frac{2ra\hat{K}(a+\rho)I}{((a+\rho)^2I^2+\hat{K})^2}+ad_h+\frac{a\beta_h}{(1+a)^2}-\frac{2r\rho \hat{K}(a+\rho)I}{((a+\rho)^2I^2+\hat{K})^2}+\frac{a^2\beta_h}{(1+a)^2}\right]\\\\
&=&\frac{\beta_h}{(1+a)^2}\left[ -\frac{2r\hat{K}(a+\rho)^2I}{((a+\rho)^2I^2+\hat{K})^2}+a d_h+\frac{a\beta_h}{1+a}\right]=\frac{\beta_h}{(1+a)^2}\left[ -\frac{2r\hat{K}(a+\rho)^2I}{((a+\rho)^2I^2+\hat{K})^2}+\tilde{d}\right].
\end{array}
\eae
If $\mathcal R_0^V=\frac{\beta_h}{d_h+\mu_h}>1$ and $\frac{r}{\tilde{d}}>\frac{2\sqrt{\hat{K}}}{a+\rho}$, then we have two interior equilibria $(S_h^k,I_h^k), k=1,2$ where
$$I_h^{1}=\frac{\frac{r}{\tilde{d}}-\sqrt{\left(\frac{r}{\tilde{d}}\right)^2-4\frac{\hat{K}}{(a+\rho)^2}}}{2},\,\,I_h^{2}=\frac{\frac{r}{\tilde{d}}+\sqrt{\left(\frac{r}{\tilde{d}}\right)^2-4\frac{\hat{K}}{(a+\rho)^2}}}{2},\,\, S_h^k= a I_h^k, k=1,2.$$
This implies that $\lambda_1(I_h^{2})\lambda_2(I_h^{1})<0,\,\, \lambda_1(I_h^{2})\lambda_2(I_h^{2})>0$ and 
$$
\begin{array}{lcl}\lambda_1(I_h^{2})+\lambda_2(I_h^{2})&<&\left(d_h+\frac{\beta_h}{1+a}\right)\left[\frac{2\hat{K}}{(a+\rho)^2I^2+\hat{K}} -1\right]\\\\
&<&\left(d_h+\frac{\beta_h}{1+a}\right)\left[\frac{2\hat{K}}{(a+\rho)^2 \frac{\hat{K}}{(a+\rho)^2}+\hat{K}} -1\right]=0\end{array}$$ since
$$I_h^{2}=\frac{\frac{r}{\tilde{d}}+\sqrt{\left(\frac{r}{\tilde{d}}\right)^2-4\frac{\hat{K}}{(a+\rho)^2}}}{2}>\frac{r}{2\tilde{d}}>\frac{\sqrt{\hat{K}}}{(a+\rho)}\Rightarrow (I_h^{2})^2>\frac{{\hat{K}}}{(a+\rho)^2}.$$
Therefore, we can conclude that if $\mathcal R_0^V=\frac{\beta_h}{d_h+\mu_h}>1$ and $\frac{r}{\tilde{d}}>\frac{2\sqrt{\hat{K}}}{a+\rho}$, then the system has two interior equilibria $(S_h^k,I_h^k), k=1,2$ where $(S_h^1,I_h^1)$ is always a saddle while $(S_h^2,I_h^2)$ is always locally asymptotically stable.\\

Notice that $\frac{r}{\tilde{d}}>\frac{2\sqrt{\hat{K}}}{a+\rho}\Leftrightarrow \frac{r}{2\sqrt{\hat{K}}}>\frac{\tilde{d}}{a+\rho}=\frac{d_h (a+1)}{a+\rho} +\frac{\mu_h}{a+\rho}>d_h$. The discussions above show that the mite-free subsystem \eqref{HV} has no interior equilibrium if either $\mathcal R_0^V<1$ or  $\frac{r}{2\sqrt{\hat{K}}}<\frac{\tilde{d}}{a+\rho}$ holds. This leads to the following two cases:
\begin{enumerate}
\item If either $\frac{r}{2\sqrt{\hat{K}}}<d_h$, then the mite-free subsystem \eqref{HV} has only the extinction equilibrium $(0,0)$ which is locally stable, thus it is globally stable.
\item If $\mathcal R_0^V>1 \mbox{ and } \,d_h< \frac{r}{2\sqrt{\hat{K}}}<\frac{\tilde{d}}{a+\rho}$, then the mite-free subsystem \eqref{HV} has the following three boundary equilibria
$$(0,0),\,\,(\bar{N}^c_h,0),\,\, \mbox{ and } (\bar{N}^*_h,0)$$ where $(0,0)$ is locally asymptotically stable; $(\bar{N}^c_h,0)$ is a source; and $(\bar{N}^*_h,0)$ is a saddle. Then according to the Poincare-Bendison theorem, we can conclude that all interior points of $\mathbb R^2_+$ converges to the extinction equilibrium $(0,0)$.
\item If $\mathcal R_0^V<1 \mbox{ and } \,d_h< \frac{r}{2\sqrt{\hat{K}}}$, then the mite-free subsystem \eqref{HV} has the following three boundary equilibria
$$(0,0),\,\,(\bar{N}^c_h,0),\,\, \mbox{ and } (\bar{N}^*_h,0)$$ where $(0,0)$ is locally asymptotically stable; $(\bar{N}^c_h,0)$ is a saddle; and $(\bar{N}^*_h,0)$ is a sink. This implies that both   $(0,0),\,\, \mbox{ and } (\bar{N}^*_h,0)$ are locally asymptotically stable.
\end{enumerate}

\end{proof}

\subsection*{Proof of Theorem \ref{th4:HVMI}}
\begin{proof}
Let $N_h=S_h+I_h$, then from Model \eqref{HVMI}, we obtain
$$\begin{array}{lcl}
N_h'&=&\frac{r\left(S_h+\rho I_h\right)^2}{\hat{K}+\left(S_h+\rho I_h\right)^2}-d_h N_h-\mu_h I_h\\\\
&\leq&\frac{rN_h^2}{\hat{K}+N_h^2}-d_h N_h
\end{array}$$ which implies that when $\frac{r}{2\sqrt{\hat{K}}}>d_h$, we have
$$\limsup_{t\rightarrow\infty}(S_h(t)+I_h(t))=\limsup_{t\rightarrow\infty}N_h(t)\leq \frac{\frac{r}{d_h}+\sqrt{\left(\frac{r}{d_h}\right)^2-4\hat{K}}}{2}=\bar{N}_h^*.$$
Therefore,  we can conclude that if the inequalities $\frac{r}{d_h}>2\sqrt{\hat{K}}$ and 
$\bar{N}_h^*<d_m+\mu_m$ hold, then we have $\limsup_{t\rightarrow\infty}S_m(t)= 0$ since
$$I_m'=c\alpha I_m\left[I_h +S_h-\frac{(d_m +\mu_m)}{c\alpha}\right]=c\alpha I_m\left[N_h-\frac{(d_m +\mu_m)}{c\alpha}\right]\leq c\alpha I_m\left[\bar{N}_h^*-\frac{(d_m +\mu_m)}{c\alpha}\right]<0.$$
This implies that the limiting dynamics of Model \eqref{HVMI} is reduced to the mite-free model \eqref{HV}.\\\\

An interior equilibrium $(S_h, I_h,I_m)$ of the healthy-mite-free subsystem \eqref{HVMI}, satisfies the following equations:
\bae\label{HVMI-eq1}
\begin{array}{lcl}
c\alpha I_m\left[I_h +S_h-\frac{d_m+\mu_m}{c\alpha}\right]&=&0\Rightarrow N_h=I_h +S_h=\frac{d_m+\mu_m}{c\alpha}\\\\
\frac{r\left(S_h+\rho I_h\right)^2}{\hat{K}+\left(S_h+\rho I_h\right)^2}-d_h N_h-\alpha N_hI_m-\mu_h I_h&=&0\Rightarrow \frac{r\left(S_h+\rho I_h\right)^2}{\hat{K}+\left(S_h+\rho I_h\right)^2}-d_h N_h-\alpha N_hI_m-\mu_h I_h=0 \\\\
&\Rightarrow&\frac{r\left(N_h-I_h+\rho I_h\right)^2}{\hat{K}+\left(N_h-I_h+\rho I_h\right)^2}-d_h N_h-\alpha N_hI_m-\mu_h I_h=0\\\\
&\Rightarrow&I_m=\frac{\frac{r\left(\frac{d_m+\mu_m}{c\alpha}-I_h+\rho I_h\right)^2}{\hat{K}+\left(\frac{d_m+\mu_m}{c\alpha}-I_h+\rho I_h\right)^2}-d_h \frac{d_m+\mu_m}{c\alpha}-\mu_h I_h}{\alpha(\frac{d_m+\mu_m}{c\alpha})}=f_1(I_h)
\end{array}
\eae
\bae\label{HVMI-eq2}
\begin{array}{lcl}
S_h\left[\frac{\beta_{h}I_h}{S_h+I_h}+\modifyr{\frac{\hat{\beta}_{mh}I_m}{S_h+I_h}}+\tilde{\beta}_{mh}I_m\right]&=&\alpha I_h I_m+(d_h+\mu_h)I_h\\\\
\left(\frac{d_m+\mu_m}{c\alpha}-I_h\right)\left[\frac{\beta_{h}I_h}{d_m +\mu_m}+\modifyr{\frac{c\alpha\hat{\beta}_{mh}I_m}{d_m+\mu_m}}+\tilde{\beta}_{mh}I_m\right]&=&\alpha I_h I_m+(d_h+\mu_h)I_h\\\\
&\Rightarrow&I_m=\frac{\left(\frac{d_m+\mu_m}{c\alpha}-I_h\right)\modifyr{\left[\frac{\beta_{h}}{d_m +\mu_m}-(d_h+\mu_h)\right]I_h}}{\alpha I_h-\left(\frac{d_m+\mu_m}{c\alpha}-I_h\right)\left(\modifyr{\frac{c\alpha\hat{\beta}_{mh}}{d_m+\mu_m}}+\tilde{\beta}_{mh}\right)}\\\\
&\Rightarrow&I_m=\frac{\left(\frac{d_m+\mu_m}{c\alpha}-I_h\right)\modifyr{\left[\frac{\beta_{h}}{d_m +\mu_m}-(d_h+\mu_h)\right]I_h}}{\left[\alpha+\modifyr{\frac{c\alpha\hat{\beta}_{mh}}{d_m+\mu_m}}+\tilde{\beta}_{mh}\right] I_h-\frac{d_m+\mu_m}{c\alpha}\left(\modifyr{\frac{c\alpha\hat{\beta}_{mh}}{d_m+\mu_m}}+\tilde{\beta}_{mh}\right)}=f_2(I_h)
\end{array}.
\eae
The equations above imply that the interior equilibrium of \eqref{HVMI} is the positive intercept of $f_1(I_h)$ and $f_2(I_h)$ subject to $0<I_h<\frac{d_m+\mu_m}{c\alpha}$. The expression for the function $f_1(I_h)$ implies that the subsystem \eqref{HVMI} has no interior equilibrium if 
$$r<\frac{d_h (d_m+\mu_m)}{c\alpha}\Leftrightarrow rc\alpha<d_h (d_m+\mu_m)$$
since $I_m=f_1(I_h)<0$ when this inequality holds. \\

Assume that $\frac{r}{2\sqrt{\hat{K}}}>\frac{d_h+\mu_h+\alpha N^*}{\rho}\mbox{ and }N_h(0)>\underline{N}^c_h$, then according to Theorem \ref{th1:basics} and \ref{th2:per-ext}, then we have 
$$\underline{N}^c_h<\underline{N}^*_h\leq \liminf_{t\rightarrow\infty}N_h(t)\leq \limsup_{t\rightarrow\infty}N_h(t)\leq \bar{N}^*_h$$ which implies that the set $\{S_h+I_h\geq \underline{N}^c_h\}$ is invariant. If $I_m=0$, then the subsystem \eqref{HVMI} reduces to the mite-free system \eqref{HV}. According to Theorem \ref{th3:MF}, the omega limit set of the mite-free system \eqref{HV} is $(\bar{N}^*_h,0)$ when $\mathcal R_0^V=\frac{\beta_h}{d_h+\mu_h}<1$ holds. If  $\bar{N}^*_h>\frac{d_m+\mu_m}{c\alpha}$, then we have
$$\begin{array}{lcl}
\frac{I_m'}{I_m}\big\vert_{I_m=0}=c\alpha\left(\bar{N}^*_h-\frac{d_m+\mu_m}{c\alpha}\right)>0.
\end{array}$$
This implies that the disease $I_m$ persists by applying the average Lynapunov theorem (Hutson 1984). Notice that 
$$\begin{array}{lcl}
I_h\big\vert_{I_h=0}=\modifyr{\frac{\hat{\beta}_{mh}S_h I_m}{N_h}}>0.
\end{array}$$
Therefore, the virus infected honeybees $I_h$ also persists.\\
\end{proof}

\subsection*{Proof of Theorem \ref{th6:global}}
\begin{proof}
The first part of Theorem \ref{th6:global} can be deduced directly from Theorem \ref{th1:basics}, \ref{th2:per-ext}, \ref{th5:DF}, \ref{th3:MF}, and \ref{th4:HVMI}. We focus on the sufficient conditions that lead to no interior equilibrium  (Item 3) and the persistence of mites (Item 4).\\

If $(S_h, I_h,S_m,I_m)$ is an interior equilibrium  of the system \eqref{newHBD}, then it satisfies the following equations:
\bae\label{newHBD-Sh}
\begin{array}{lcl}
0&=&S_m\left[c\alpha S_h-\frac{\hat{\beta}_{hm}I_h}{S_h+I_h}-d_m\right]\\\\
&\Rightarrow&I_h=S_h \frac{c\alpha S_h-d_m}{d_m+\hat{\beta}_{hm}-c\alpha S_h},  \frac{d_m}{c\alpha} <S_h<\frac{d_m+\hat{\beta}_{hm}}{c\alpha} \\\\
&\Rightarrow&N_h=S_h+I_h=\frac{\hat{\beta}_{hm} S_h}{d_m+\hat{\beta}_{hm}-c\alpha S_h}=g_1(S_h)
\end{array}
\eae
\bae\label{newHBD-Ih}
\begin{array}{lcl}
0&=& \frac{r\left(S_h+\rho I_h\right)^2}{\hat{K}+\left(S_h+\rho I_h\right)^2}-d_h (S_h+I_h)-\alpha (S_h+I_h) (S_m+I_m)-\mu_h I_h\\\\
&\Rightarrow&N_m=S_m+I_m=\frac{\frac{r\left(S_h+\rho I_h\right)^2}{\hat{K}+\left(S_h+\rho I_h\right)^2}-d_h (S_h+I_h)-\mu_h I_h}{\alpha (S_h+I_h)}\\\\
&\Rightarrow&N_m=\frac{rS_h\left(\hat{\beta}_{hm}+(d_m-c\alpha S_h)\right)\left(\hat{\beta}_{hm}+(d_m-c\alpha S_h)(1-\rho)\right)^2}{\alpha\hat{\beta}_{hm}\left[\hat{K}\left(\hat{\beta}_{hm}+(d_m-c\alpha S_h)\right)^2+S_h^2\left(\hat{\beta}_{hm}+(d_m-c\alpha S_h)(1-\rho)\right)^2\right]}-\frac{d_h}{\alpha} -\frac{\mu_h(c\alpha S_h-d_m)}{\hat{\beta}_{hm}} =g_2(S_h)
\end{array}
\eae\bae\label{newHBD-Sm}
\begin{array}{lcl}
0&=&c\alpha(S_h+I_h) (S_m+I_m)-d_m(S_m+I_m) -\mu_mI_m\\\\
&\Rightarrow& I_m=\frac{N_m(c\alpha(S_h+I_h) -d_m)}{\mu_m}=\frac{N_m\left[c\alpha\hat{\beta}_{hm} S_h -d_m\left(\hat{\beta}_{hm}+(d_m-c\alpha S_h)\right)\right]}{\mu_m\left(\hat{\beta}_{hm}+(d_m-c\alpha S_h)\right)}=\frac{N_m(c\alpha S_h-d_m)(\hat{\beta}_{hm}+d_m) }{\mu_m\left(\hat{\beta}_{hm}+(d_m-c\alpha S_h)\right)}\\\\
&&=N_mg_3(S_h)=g_2(S_h)g_3(S_h)
\end{array}
\eae

{\small\bae\label{newHBD-Im}
\begin{array}{lcl}
0&=&S_h\left[\frac{\beta_{h}I_h}{S_h+I_h}+\modifyr{\frac{\hat{\beta}_{mh}I_m}{S_h+I_h}}+\tilde{\beta}_{mh}I_m\right]-\alpha I_h (S_m+I_m)-(d_h+\mu_h)I_h\\\\
&\Rightarrow&\frac{d_m+\hat{\beta}_{hm}-c\alpha S_h}{c\alpha S_h-d_m}\left[\frac{\beta_{h}(c\alpha S_h-d_m)}{\hat{\beta}_{hm}}+\modifyr{\hat{\beta}_{mh}\frac{N_m(c\alpha S_h-d_m)(\hat{\beta}_{hm}+d_m) }{\mu_m \hat{\beta}_{hm}S_h}}+\tilde{\beta}_{mh}\frac{N_m(c\alpha S_h-d_m)(\hat{\beta}_{hm}+d_m) }{\mu_m\left(\hat{\beta}_{hm}+(d_m-c\alpha S_h)\right)}\right]-\alpha N_m-(d_h+\mu_h)=0\\\\
&\Rightarrow&\frac{\beta_{h}(d_m+\hat{\beta}_{hm}-c\alpha S_h)}{\hat{\beta}_{hm}}+\modifyr{\frac{\hat{\beta}_{mh}N_m(\hat{\beta}_{hm}+d_m)(d_m+\hat{\beta}_{hm}-c\alpha S_h) }{\mu_m \hat{\beta}_{hm}S_h}}+\frac{\tilde{\beta}_{mh}N_m(\hat{\beta}_{hm}+d_m) }{\mu_m}-(d_h+\mu_h)=\alpha N_m\\\\
&\Rightarrow&N_m=\frac{\frac{\beta_{h}(d_m+\hat{\beta}_{hm}-c\alpha S_h)}{\hat{\beta}_{hm}}-(d_h+\mu_h)}{\alpha-\modifyr{\frac{\hat{\beta}_{mh}(\hat{\beta}_{hm}+d_m)(d_m+\hat{\beta}_{hm}-c\alpha S_h) }{\mu_m \hat{\beta}_{hm}S_h}}-\frac{\tilde{\beta}_{mh}(\hat{\beta}_{hm}+d_m) }{\mu_m}}=g_4(S_h)
\end{array}.
\eae}
Therefore, the interior equilibrium $(S_h, I_h,S_m,I_m)$ are positive solutions of the following four equations:
{\small\bae\label{newHBD-Im2}
\begin{array}{lcl}
N_h&=&S_h+I_h=g_1(S_h)=\frac{\hat{\beta}_{hm} S_h}{d_m+\hat{\beta}_{hm}-c\alpha S_h}\mbox{ subject to } S_h<\frac{d_m+\hat{\beta}_{hm}}{c\alpha}\\\\
\frac{I_m}{N_m}&=&g_3(S_h)=\frac{(c\alpha S_h-d_m)(\hat{\beta}_{hm}+d_m) }{\mu_m\left(\hat{\beta}_{hm}+d_m-c\alpha S_h\right)}<1 \mbox{ subject to }\frac{d_m}{c\alpha}<S_h<\frac{d_m+\hat{\beta}_{hm}}{c\alpha}\frac{d_m+\mu_m}{d_m+\hat{\beta}_{hm}+\mu_m}\\\\
N_m&=&g_2(S_h)=\frac{rS_h\left(\hat{\beta}_{hm}+(d_m-c\alpha S_h)\right)\left(\hat{\beta}_{hm}+(d_m-c\alpha S_h)(1-\rho)\right)^2}{\alpha\hat{\beta}_{hm}\left[\hat{K}\left(\hat{\beta}_{hm}+(d_m-c\alpha S_h)\right)^2+S_h^2\left(\hat{\beta}_{hm}+(d_m-c\alpha S_h)(1-\rho)\right)^2\right]}-\frac{d_h}{\alpha} -\frac{\mu_h(c\alpha S_h-d_m)}{\hat{\beta}_{hm}} \\\\
N_m&=&g_4(S_h)=\frac{\frac{\beta_{h}(d_m+\hat{\beta}_{hm}-c\alpha S_h)}{\hat{\beta}_{hm}}-(d_h+\mu_h)}{\alpha-\modifyr{\frac{\hat{\beta}_{mh}(\hat{\beta}_{hm}+d_m)(d_m+\hat{\beta}_{hm}-c\alpha S_h) }{\mu_m \hat{\beta}_{hm}S_h}}-\frac{\tilde{\beta}_{mh}(\hat{\beta}_{hm}+d_m) }{\mu_m}}\\\\
&=&
\frac{\mu_m S_h\left[\beta_{h}(d_m+\hat{\beta}_{hm})-\hat{\beta}_{hm}(d_h+\mu_h)-c\alpha \beta_{h}S_h\right]}{\left[\alpha\mu_m \hat{\beta}_{hm}+c\alpha \hat{\beta}_{mh}(\hat{\beta}_{hm}+d_m)-\tilde{\beta}_{mh} \hat{\beta}_{hm}(\hat{\beta}_{hm}+d_m)\right]S_h-\modifyr{\hat{\beta}_{mh}(\hat{\beta}_{hm}+d_m)(d_m+\hat{\beta}_{hm})}}
\end{array}.
\eae}
Thus, we can conclude the following statements regarding the sign of $g_4(S_h)$:
\begin{enumerate}
\item If $\beta_{h}(d_m+\hat{\beta}_{hm})<\hat{\beta}_{hm}(d_h+\mu_h)$ and $\alpha\mu_m \hat{\beta}_{hm}+c\alpha \hat{\beta}_{mh}(\hat{\beta}_{hm}+d_m)<\tilde{\beta}_{mh} \hat{\beta}_{hm}(\hat{\beta}_{hm}+d_m)$, then $g_4(S_h)>0$ for all $S_h>0$.
\item If $\beta_{h}(d_m+\hat{\beta}_{hm})>\hat{\beta}_{hm}(d_h+\mu_h)$ and $\alpha\mu_m \hat{\beta}_{hm}+c\alpha \hat{\beta}_{mh}(\hat{\beta}_{hm}+d_m)>\tilde{\beta}_{mh} \hat{\beta}_{hm}(\hat{\beta}_{hm}+d_m)$, then $g_4(S_h)>0$ when 
$$\frac{\hat{\beta}_{mh}(\hat{\beta}_{hm}+d_m)(d_m+\hat{\beta}_{hm})}{\alpha\mu_m \hat{\beta}_{hm}+c\alpha \hat{\beta}_{mh}(\hat{\beta}_{hm}+d_m)-\tilde{\beta}_{mh} \hat{\beta}_{hm}(\hat{\beta}_{hm}+d_m)}<S_h<\frac{\beta_{h}(d_m+\hat{\beta}_{hm})-\hat{\beta}_{hm}(d_h+\mu_h)}{c\alpha \beta_{h}}.$$
\item If $\beta_{h}(d_m+\hat{\beta}_{hm})<\hat{\beta}_{hm}(d_h+\mu_h)$ and $\alpha\mu_m \hat{\beta}_{hm}+c\alpha \hat{\beta}_{mh}(\hat{\beta}_{hm}+d_m)>\tilde{\beta}_{mh} \hat{\beta}_{hm}(\hat{\beta}_{hm}+d_m)$, then $g_4(S_h)>0$ when 
$$S_h<\frac{\hat{\beta}_{mh}(\hat{\beta}_{hm}+d_m)(d_m+\hat{\beta}_{hm})}{\alpha\mu_m \hat{\beta}_{hm}+c\alpha \hat{\beta}_{mh}(\hat{\beta}_{hm}+d_m)-\tilde{\beta}_{mh} \hat{\beta}_{hm}(\hat{\beta}_{hm}+d_m)}.$$
\item If $\beta_{h}(d_m+\hat{\beta}_{hm})>\hat{\beta}_{hm}(d_h+\mu_h)$ and $\alpha\mu_m \hat{\beta}_{hm}+c\alpha \hat{\beta}_{mh}(\hat{\beta}_{hm}+d_m)<\tilde{\beta}_{mh} \hat{\beta}_{hm}(\hat{\beta}_{hm}+d_m)$, then $g_4(S_h)>0$ when $$S_h>\frac{\beta_{h}(d_m+\hat{\beta}_{hm})-\hat{\beta}_{hm}(d_h+\mu_h)}{c\alpha \beta_{h}}.$$
\end{enumerate}
Notice that the interior equilibrium $(S_h, I_h,S_m,I_m)$ requires $g_3(S_h)>0$, i.e., $\frac{d_m}{c\alpha}<S_h<\frac{d_m+\hat{\beta}_{hm}}{c\alpha}\frac{d_m+\mu_m}{d_m+\hat{\beta}_{hm}+\mu_m}$. Therefore, the interior equilibrium $(S_h, I_h,S_m,I_m)$ does not exist if one of the following inequalities hold
\begin{enumerate}
\item $\frac{\beta_{h}}{d_h+\mu_h}>\frac{\hat{\beta}_{hm}}{d_m+\hat{\beta}_{hm}},\,\,\frac{\alpha\mu_m \hat{\beta}_{hm}}{\hat{\beta}_{hm}+d_m}> \tilde{\beta}_{mh} \hat{\beta}_{hm}-c\alpha$ and
$$\frac{\hat{\beta}_{mh}}{\frac{\alpha\mu_m \hat{\beta}_{hm}}{(d_m+\hat{\beta}_{hm})}+c\alpha \hat{\beta}_{mh}-\tilde{\beta}_{mh} \hat{\beta}_{hm}}>\frac{d_m+\mu_m}{c\alpha(d_m+\hat{\beta}_{hm}+\mu_m)}\mbox{ or } \frac{\frac{\beta_{h}}{d_h+\mu_h}-\frac{\hat{\beta}_{hm}}{d_m+\hat{\beta}_{hm}}}{ \beta_{h}(d_m+\hat{\beta}_{hm})(d_h+\mu_h)}<d_m.$$ 
\item 
$\frac{\beta_{h}}{d_h+\mu_h}<\frac{\hat{\beta}_{hm}}{d_m+\hat{\beta}_{hm}},\,\,\frac{\alpha\mu_m \hat{\beta}_{hm}}{\hat{\beta}_{hm}+d_m}> \tilde{\beta}_{mh} \hat{\beta}_{hm}-c\alpha$ and
$\frac{\hat{\beta}_{mh}(d_m+\hat{\beta}_{hm})}{\frac{\alpha\mu_m \hat{\beta}_{hm}}{(d_m+\hat{\beta}_{hm})}+c\alpha \hat{\beta}_{mh}-\tilde{\beta}_{mh} \hat{\beta}_{hm}}<\frac{d_m}{c\alpha}.$ 
\item $\frac{\beta_{h}}{d_h+\mu_h}>\frac{\hat{\beta}_{hm}}{d_m+\hat{\beta}_{hm}},\,\,\frac{\alpha\mu_m \hat{\beta}_{hm}}{\hat{\beta}_{hm}+d_m}< \tilde{\beta}_{mh} \hat{\beta}_{hm}-c\alpha$ and $\frac{\frac{\beta_{h}}{d_h+\mu_h}-\frac{\hat{\beta}_{hm}}{d_m+\hat{\beta}_{hm}}}{ \beta_{h}(d_m+\hat{\beta}_{hm})(d_h+\mu_h)}>\frac{(d_m+\mu_m)(d_m+\hat{\beta}_{hm})}{d_m+\hat{\beta}_{hm}+\mu_m}.$\\
\end{enumerate}

Assume that $\frac{r}{2\sqrt{\hat{K}}}>\frac{d_h+\mu_h+\alpha N^*}{\rho}\mbox{ and }N_h(0)>\underline{N}^c_h$, then according to Theorem \ref{th1:basics} and \ref{th2:per-ext}, then we have 
$$\underline{N}^c_h<\underline{N}^*_h\leq \liminf_{t\rightarrow\infty}N_h(t)\leq \limsup_{t\rightarrow\infty}N_h(t)\leq \bar{N}^*_h$$ which implies that the set $\{S_h+I_h\geq \underline{N}^c_h\}$ is invariant. If $S_m=I_m=0$, then the  full system \eqref{newHBD} reduces to the mite-free system \eqref{HV}. According to Theorem \ref{th3:MF}, the omega limit set of the mite-free system \eqref{HV} is $(\bar{N}^*_h,0)$ when $\mathcal R_0^V=\frac{\beta_h}{d_h+\mu_h}<1$ holds. If  $\bar{N}^*_h>\frac{d_m+\mu_m}{c\alpha}$, then we have
$$\begin{array}{lcl}
\frac{N_m'}{N_m}\big\vert_{N_m=0}\geq c\alpha\left(\bar{N}^*_h-\frac{d_m+\mu_m}{c\alpha}\right)>0.
\end{array}$$
This implies that the mite population $N_m$ persists by applying the average Lynapunov theorem (Hutson 1984).\\

\end{proof}
\section*{Acknowledgements}
Y.K's research is partially supported by NSF-DMS (1313312). \modifyc{We really appreciate useful comments and suggestions provided by two reviewers to significantly improve this manuscript.}\\
\section*{References}
\bibliographystyle{plain}
\bibliography{references}

\end{document}